

\input amstex
\documentstyle{amsppt}

\input label.def
\input degt.def


\input epsf
\def\picture#1{\epsffile{#1-bb.eps}}
\def\cpic#1{$\vcenter{\hbox{\picture{#1}}}$}

\def\ie{\emph{i.e\.}}
\def\eg{\emph{e.g\.}}
\def\cf{\emph{cf\.}}
\def\via{\emph{via}}
\def\etc{\emph{etc}}

{\catcode`\@11
\gdef\proclaimfont@{\sl}
\gdef\subsubheadfont@{\subheadfont@}
}

\newtheorem{Remark}\Remark\thm\endAmSdef
\conjecture\thm\endproclaim
\question\thm\endAmSdef
\problem\thm\endAmSdef
\example\thm\endAmSdef
\newhead\subsubsection\subsubsection\endsubhead
\def\paragraph{\subsubsection{}}

\def\dash{\item"\hfill--\hfill"}
\def\Dashes{\widestnumber\item{--}\roster}
\def\endDashes{\endroster}

\loadbold
\def\bA{\bold A}
\def\bD{\bold D}
\def\bE{\bold E}

\def\I{{\fam0I}}
\def\II{{\fam0II}}
\def\III{{\fam0III}}
\def\IV{{\fam0IV}}
\def\nii{n_{\II}}
\def\niii{n_{\III}}
\def\niv{n_{\IV}}

\let\Ga\alpha
\let\Gb\beta
\let\Gg\gamma

\let\Gf\varphi
\let\Gs\sigma

\let\Gr\rho
\def\1{^{-1}}
\def\bv{\bold v}

\let\splus\oplus

\def\bT{T}
\def\tT{\topsmash{\tilde T}}
\let\<\langle
\let\>\rangle
\let\onto\twoheadrightarrow
\let\into\hookrightarrow
\let\gtimes\otimes
\def\ls|#1|{\mathopen|#1\mathclose|}

\def\GL{\operatorname{\text{\sl GL\/}}}
\def\SL{\operatorname{\text{\sl SL\/}}}
\def\PGL{\operatorname{\text{\sl PGL\/}}}
\def\PSL{\operatorname{\text{\sl PSL\/}}}
\def\NS{\operatorname{\text{\sl NS\/}}}
\def\MW{\operatorname{\text{\sl MW\/}}}
\def\Sp{\operatorname{Sp}}
\def\BG#1{\Bbb B_{#1}}
\def\BGM#1{\frak B_{#1}}
\def\bmtimes{\mathbin{\bold\cdot}}
\def\Aut{\operatorname{Aut}}

\def\gm{\frak m}
\def\gM{\bar\gm}
\def\BM-{BM-\hskip0pt }
\def\bminf{\gm_\infty}
\def\bmP#1{\bminf(#1)}    
\def\bmH#1{\operatorname{Im}(#1)}      
\def\bmTL{\Cal T}      
\def\bmT#1{\bmTL(#1)}  
\def\bmpi#1{\pi_1(#1)} 
\let\MG\Gamma          
\def\bMG{\topsmash{\tilde\Gamma}}   
\def\bH{\topsmash{\tilde H}}        
\def\SG{\Bbb S}        
\def\CG#1{\Z_{#1}}     
\def\Cp#1{\Bbb P^{#1}}
\def\Rp#1{\Cp{#1}_\R}

\def\Sk{\operatorname{Sk}}
\def\Ski{\Sk^\circ}
\def\Skii{\Sk^\bullet}
\def\Skc{\Sk^c}
\def\Stab{\operatorname{Stab}}
\def\op{\operatorname{op}}
\def\nx{\operatorname{nx}}
\def\gO{\frak O}
\def\val{\operatorname{val}}
\def\X{\Bbb X}
\def\Y{\Bbb Y}
\def\CH{\Cal H}
\def\CL{\Cal L}
\def\CM{\Cal M}
\def\CV{\Cal V}
\def\CP{\Cal P}
\def\CE{\Cal E}
\def\CESk{\CE_{\Sk}}
\def\ee{{\fam0e}}

\def\cc[#1]{[\![#1]\!]}  

\let\XX=X                      
\def\Xs{\XX^\sharp}            
\let\BB=B                      
\def\Bs{\BB^\sharp}            
\def\jX{j_{\XX}}               
\def\hh{\frak h}               
\def\bhh{\tilde\frak h}        
\def\hX{\hh_{\XX}}             
\def\bhX{\topsmash{\tilde\hX}} 
\def\SkX{\Sk_{\XX}}            
\let\CC=C                      
\def\jC{j_{\CC}}               
\def\CCa{\CC^{\fam0a}}         
\let\bm\mu                     
\def\bmC{\bm_{\CC}}            
\def\gMC{\gM_{\CC}}

\def\barF{\bar F}
\def\piorb{\pi^{\fam0orb}_1}
\let\tree\Xi

\def\GAP{{\tt GAP}}

\def\inserthyphen{\ifcat\next a-\fi\ignorespaces}
\let\BLACK\bullet
\let\WHITE\circ
\def\CROSS{\vcenter{\hbox{$\scriptstyle\mathord\times$}}}
\let\STAR*
\def\TRIANG{\vcenter{\hbox{$\scriptstyle\mathord\vartriangle$}}}
\def\pblack-{$\BLACK$\futurelet\next\inserthyphen}
\def\pwhite-{$\WHITE$\futurelet\next\inserthyphen}
\def\pcross-{$\CROSS$\futurelet\next\inserthyphen}
\def\pstar-{$\STAR$\futurelet\next\inserthyphen}
\def\ptriang-{$\TRIANG$\futurelet\next\inserthyphen}
\def\black{\protect\pblack}
\def\white{\protect\pwhite}
\def\cross{\protect\pcross}
\def\star{\protect\pstar}
\def\triang{\protect\ptriang}
\def\NO#1{n_{#1}}
\def\nblack{\NO\BLACK}
\def\nwhite{\NO\WHITE}

\def\nstar{\NO\STAR}
\def\ntriang{\NO\TRIANG}

\def\MGSk#1{{\BLACK}{\joinrel\relbar\joinrel\relbar\joinrel}{#1}}
\def\stone#1{$\joinrel\relbar\joinrel\mathrel{\vcenter{\hbox{$#1$}}}\joinrel\relbar\joinrel$}
\def\stone#1{${\relbar\joinrel}{\vcenter{\hbox{$#1$}}}{\joinrel\relbar}$}
\def\stcirc{\scriptstyle\bigcirc}
\def\stleft{\rlap{$<$}{\relbar}}
\def\stright{\rlap{$>$}{\relbar}}

\topmatter

\author
Alex Degtyarev
\endauthor

\title
Hurwitz equivalence of braid monodromies and extremal elliptic surfaces
\endtitle

\rightheadtext{Hurwitz equivalence of braid monodromies}

\address
Department of Mathematics,
Bilkent University,
06800 Ankara, Turkey
\endaddress

\email
degt\@fen.bilkent.edu.tr
\endemail

\abstract
We discuss the equivalence between the categories of
certain ribbon graphs and subgroups of the
modular group~$\Gamma$ and use it to construct exponentially large
families of not Hurwitz equivalent simple braid monodromy
factorizations
of the same element.
As an application, we also obtain exponentially large families
of {\it topologically\/} distinct algebraic objects such as
extremal elliptic surfaces, real trigonal curves, and real
elliptic surfaces.
\endabstract



\endtopmatter

\document

\section{Introduction\label{S.intro}}

Strictly speaking, principal results of the paper concern extremal
elliptic surfaces, see Subsection~\ref{s.intro.elliptic}. However,
we start with discussing a few applications to
the braid monodromy,
which seems to be a subject of more general interest.

\subsection{Braid monodromy}\label{s.intro.BM}
Throughout the paper, we use the notation
$\cc[\,\cdot\,]=\cc[\,\cdot\,]_G$ for the conjugacy class of an
element $g\in G$ or a subgroup $H\subset G$ of a group~$G$.

\definition\label{def.BMF}
Given a group~$G$, a
\emph{\rom($G$-valued\rom) braid monodromy factorization}
(\emph{\BM-factorization} for short) of length~$r$ is a
sequence $\gM=(\gm_1,\ldots,\gm_r)$ of elements of~$G$. Two
\BM-factorizations are \emph{strongly \rom(Hurwitz\rom)
equivalent} if they are related by a finite sequence of
\emph{Hurwitz moves}
$$
(\ldots,\gm_i,\gm_{i+1},\ldots)\mapsto
 (\ldots,\gm_i\1\gm_{i+1}\gm_i,\gm_i,\ldots)
$$
and their inverse. Two \BM-factorizations are \emph{weakly
equivalent} if they are related by a sequence of Hurwitz moves and
their inverse
and/or \emph{global conjugation}
$$
\gM=(\gm_i)\mapsto g\1\gM g:=(g\1\gm_ig),\quad g\in G.
$$

Often it is required that each element~$\gm_i$ of a
\BM-factorization should belong to the union $\bigcup_jC_j$ of
several conjugacy classes~$C_j$ fixed in advance. Thus, a
$\BG{n}$-valued \BM-factorization is called \emph{simple} if
each~$\gm_i$ is conjugate to the Artin generator
$\Gs_1$,
see Definition~\ref{def.simple}.
\enddefinition

Note that we regard a braid monodromy as an
\emph{anti}-homomorphism, see~\ref{s.geometry} below. This
convention explains the slightly unusual form of the Hurwitz moves
and the fact that the order of multiplication is reversed
in~\iref{s.invariants}{bmP}.

In this paper we
mainly
deal with
the first nonabelian
braid group~$\BG3$ and the closely related groups
$\bMG:=\SL(2,\Z)$ and $\MG:=\PSL(2,\Z)$. A $\bMG$- or $\MG$-valued
\BM-factorization $(\gm_i)$ is called \emph{simple} if
each~$\gm_i$ belongs to the
corresponding
conjugacy class $\cc[\X\Y]$, see
Subsection~\ref{s.MG} for the notation. The classifications of
simple \BM-factorizations (up to weak/strong Hurwitz equivalence)
in all
three groups
coincide,
see Proposition~\ref{BG3=MG}.

\paragraph\label{s.geometry}
A $G$-valued \BM-factorization $\gM=(\gm_1,\ldots,\gm_r)$ can be
regarded as an anti-homo\-morphism
$\<\Gg_1,\ldots,\Gg_r\>\to G$, $\Gg_i\mapsto\gm_i$,
$i=1,\ldots,r$. In this interpretation,
Hurwitz moves generate the canonical action of the braid
group~$\BG{r}$ on the free group $\<\Gg_1,\ldots,\Gg_r\>$, and the
global conjugation represents the adjoint action of~$G$ on itself.
Geometrically, anti-homomorphisms as above arise from locally
trivial fibrations $\Xs\to\Bs$ over a punctured disk; then $G$ is
the (appropriately defined) mapping class group of the fiber over
a fixed point $b\in\partial\Bs$
and
$\<\Gg_1,\ldots,\Gg_r\>$ is a geometric basis for
$\pi_1(\Bs,b)$. In this set-up, Hurwitz moves can be
interpreted either as basis changes or as automorphisms of~$\Bs$
fixed on the boundary, see~\cite{Artin}, and the topological
classification of fibrations reduces to the purely algebraic
classification of $G$-valued \BM-factorizations up to weak Hurwitz
equivalence.
The best known examples are
\Dashes
\dash
ramified coverings (the fiber is a finite set and
$G=\SG_n$, see~\cite{Hurwitz});
\dash
algebraic or, more generally,
pseudoholomorphic and Hurwitz curves in~$\C^2$
(the fiber is a punctured plane and
$G=\BG{n}$, see~\cite{Zariski.group}, \cite{vanKampen},
\cite{Chisini1}, \cite{Chisini2},
\cite{Moishezon1}, \cite{Moishezon2},
\cite{Kulikov}, \cite{KK},
\cite{Orevkov}, \cite{Orevkov.talk});
\dash
(real) elliptic surfaces or, more generally,
(real) Lefschetz fibrations of genus one
(the fiber is
an elliptic curve/topological torus
and $G=\bMG$, see~\cite{Kodaira}, \cite{Shioda.modular},
\cite{Moishezon1},
\cite{Bogomolov},
\cite{dessin-e7}, \cite{DIK.elliptic},
\cite{Orevkov}, \cite{Orevkov.talk},
\cite{Nermin}).
\endDashes
Last two subjects are quite popular and the reference
lists are far from complete: I tried to cite the founding papers
and a few recent results/surveys only.

Usually it is understood that the punctures of~$\Bs$ correspond to
the singular fibers of a fibration~$\XX\to\BB$ over a disk, the
type of each singular fiber~$F$ being represented by the conjugacy
class of the local monodromy about~$F$. Thus, in the three
examples above, simple \BM-factorizations correspond to fibrations
with simplest, not removable by a small deformation, singular
fibers.

\paragraph\label{s.invariants}
The following is a list of the most commonly used
weak/strong equivalence invariants of a
$G$-valued \BM-factorization~$\gM$:
\roster
\item\local{bmP}
the \emph{monodromy at infinity} $\bmP\gM:=\gm_r\ldots\gm_1\in G$
is a strong invariant;
its conjugacy class $\cc[\bmP\gM]$
is a weak invariant;
\item\local{bmH}
the \emph{monodromy group}
$\bmH\gM:=\<\gm_1,\ldots,\gm_r\>\subset G$
is a strong invariant; its
conjugacy class $\cc[\bmH\gM]$ is a weak invariant;
\item\local{bmT}
for
$G=\SL(2,\Z)$,
the \emph{transcendental lattice} $\bmT\gM$,
see Subsection~\ref{s.lattice} for
the definition and generalizations,
is a week invariant;
\item\local{bmpi}
for $G=\BG3$, define the
\emph{\rom(affine\rom) fundamental group}
(see~\cite{Zariski.group}, \cite{vanKampen})
$$
\bmpi\gM:=\<\Ga_1,\Ga_2,\Ga_3\,|\,
 \text{$\gm_i(\Ga_j)=\Ga_j$ for $i=1,\ldots,r$, $j=1,2,3$}\>;
$$
the homomorphism $\<\Ga_1,\Ga_2,\Ga_3\>\onto\bmpi\gM$ is a strong
invariant; it depends on $\bmH\gM$ only; the isomorphism class of
the abstract group
$\bmpi\gM$ is a weak invariant; it depends on $\cc[\bmH\gM]$ only.
\endroster
Due to Proposition~\ref{BG3=MG}, invariants~\loccit{bmT}
and~\loccit{bmpi} apply equally well to simple
\hbox{$\BG3$-}, \hbox{$\bMG$-}, and $\MG$-valued
\BM-factorizations. Note that often it is the group~\loccit{bmpi}
that is the ultimate goal of computing the \BM-factorization in
the first place.

Geometrically, most important is the monodromy at
infinity~\loccit{bmP}; in the set-up of~\ref{s.geometry}, it
corresponds to the monodromy along the boundary~$\partial\BB$, and
the \BM-factorizations~$\gM$ with a given class
$\cc[\bmP\gM]\subset G$
enumerate the extensions to~$\BB$ of a given fibration
over~$\partial\BB$. For this reason, a \BM-factorization~$\gM$ is
often regarded as a factorization of a given element $\bmP\gM$
(which explains the term).
The geometric importance of the extension problem, a number of
partial results, and extensive experimental evidence give rise to
the following two long standing questions.

\question\label{?.uniqueness}
Is the weak/strong equivalence class of a simple $\BG{n}$-valued
\BM-factorization~$\gM$ determined by the monodromy at infinity
$\bmP\gM$?
(Note that the length of~$\gM$ is determined
by~$\bmP\gM$, see~\ref{s.length}.)
\endquestion

\question\label{?.weak=strong}
If two
simple $\BG{n}$-valued \BM-factorizations $\gM_1$, $\gM_2$
have the same monodromy at infinity
and are weakly equivalent, are they also
strongly equivalent? In other words, if a simple
\BM-factorization~$\gM$ is conjugated
by an element of~$G$ commuting with
$\bmP\gM$, is the result strongly equivalent to~$\gM$?
\endquestion

The answer to Question~\ref{?.uniqueness} is in the affirmative if
$n=3$ and $\bmP\gM$ is a central (see~\cite{Moishezon1})
or, more generally, positive (with respect to the Artin basis,
see~\cite{Orevkov.talk}) element
of~$\BG3$. Furthermore, for any~$n$, two \BM-factorizations
sharing the same monodromy at infinity are known to be
\emph{stably equivalent}, see~\cite{KK} or~\cite{Kulikov} for
details.

The condition that $\gM$ should be simple in
Question~\ref{?.uniqueness} is crucial: in general, a
\BM-factorization is not unique. First example was
essentially found in~\cite{Zariski.group}, and a great deal
of other examples have been discovered since then. A few new
examples are discussed in Subsections~\ref{s.more.BM}
and~\ref{s.sextics}. In particular, we give a very simple, not
computer aided, proof of the non-equivalence of the two
\BM-factorizations considered in~\cite{Artal.braids}.

\subsection{Principal results}\label{s.intro.results}
We answer
Questions~\ref{?.uniqueness} and~\ref{?.weak=strong}
in the negative for
the braid group $\BG3$ (and
related groups~$\MG$
and~$\bMG$, see Proposition~\ref{BG3=MG}). The
inclusion $\BG3\into\BG{n}$ implies a
negative answer for the other braid groups as well, at least
concerning the strong equivalence, see~\ref{s.B3toBn}.

Let $\bT(k)$ be the number of
isotopy
classes of
trees $\tree\subset S^2$ with
$k$ trivalent vertices and $(k+2)$ monovalent vertices (and no
other vertices), see Section~\ref{S.trees} and
Corollary~\ref{cor.counts}.
Let $C(k)=\binom{2k}{k}/(k+1)$ be the $k$-th Catalan number, and
let $\tT(k)=(5k+4)C(k)/(k+2)$, see Subsection~\ref{s.counts} and
Corollary~\ref{cor.counts}. Note that each of the three series
grows faster that $a^k$ for any $a<4$. The first few values
of $\bT(k)$ and $\tT(k)$ are listed in Table~\ref{tab.T}.

\midinsert
\table\label{tab.T}
A few values of $\bT(k)$ and $\tT(k)$
\endtable
\centerline{\vbox{\halign{\strut\hss$#$\hss&&\quad\hss$#$\hss\cr
k&0&1&2&3&4&5&6&7&\dots&10&\dots&15\cr
\bT(k)&1&1&1&1&4&6&19&49&\dots&1424&\dots&570285\cr
\tT(k)&2&3&7&19&56&174&561&1859&\dots&75582&\dots&45052515\crcr}}}
\endinsert

\theorem\label{th.strong}
For each integer $k\ge0$, there is a set $\{\gM_i\}$,
$i=1,\ldots,\tT(k)$,
of simple $\MG$-valued \BM-factorizations of
length $(k+2)$ that share the same
\Dashes
\dash
monodromy at infinity $\bmP{\gM_i}=(\X\Y)^{-5k-4}$,
\dash
transcendental lattice $\bmT{\gM_i}$, see
Example~\ref{ex.lattice}, and
\dash
fundamental group $\bmpi{\gM_i}$ \rom(which is~$\Z$ for
$k\ge2$\rom)
\endDashes
but are not strongly equivalent\rom:
the monodromy
groups $\bmH{\gM_i}\subset\MG$ are pairwise distinct
subgroups of index~$6(k+1)$.
\endtheorem

\theorem\label{th.weak}
For each~$k$,
the \BM-factorizations~$\gM_i$ in Theorem~\ref{th.strong} form
$\bT(k)$ distinct weak equivalence classes\rom: they
are distinguished by the
conjugacy classes $\cc[\bmH{\gM_i}]$ of the monodromy groups.
\endtheorem

Since $\bT(k)<\tT(k)$ for all $k\ge0$,
one has the following corollary.

\corollary\label{weak=/=strong}
For each integer $k\ge0$, there is a pair
of
conjugate simple $\MG$-valued \BM-factorizations of
length~$(k+2)$ that share the same monodromy at infinity
$(\X\Y)^{-5k-4}$
but are not strongly equivalent.
\qed
\endcorollary

Theorems~\ref{th.strong} and~\ref{th.weak} are proved in
Subsection~\ref{s.proof.weak}; the \BM-factorizations in question
are given by~\eqref{eq.MG}, and their $\BG3$-valued counterparts
are given by~\eqref{eq.BG}.
The first example of weakly but not strongly equivalent
$\BG3$-valued
\BM-factorizations given by Corollary~\ref{weak=/=strong} has
length two; it is as simple as
$$
\gM'=(\Gs_1^2\Gs_2\Gs_1^{-2},\Gs_2),\qquad
\gM''=(\Gs_1\Gs_2\Gs_1\1,\Gs_1\1\Gs_2\Gs_1),
$$
see Example~\ref{ex.k=0}. The first example of non-equivalent
\BM-factorizations given by Theorem~\ref{th.weak} has length six,
see Example~\ref{ex.k=4}. In Subsection~\ref{s.two} we construct
another example of not weakly equivalent \BM-factorizations of
length two; they also differ by the monodromy groups, which are of
infinite index. A few other examples
(not necessarily simple)
are considered in
Subsections~\ref{s.more.BM} and~\ref{s.sextics}.

\subsection{Elliptic surfaces}\label{s.intro.elliptic}
Recall that an \emph{extremal elliptic surface} can be defined
as a Jacobian elliptic
surface~$\XX$ of maximal Picard number,
$\rank\NS(\XX)=h^{1,1}(\XX)$,
and minimal Mordell-Weil rank, $\rank\MW(\XX)=0$.
(For an alternative description, in terms of singular fibers,
see~\ref{s.extremal}. Yet another characterization is the
following: a Jacobian elliptic surface is extremal if and only if
its transcendental lattice is positive definite,
see~\cite{tripods}.)
Extremal elliptic surfaces are rigid
(any
small
fiberwise equisingular deformation of such a surface~$\XX$ is
isomorphic to~$\XX$);
they are defined over
algebraic number fields.

In this paper, we mainly deal with elliptic surfaces with singular
fibers of Kodaira types~$\I_p$ and~$\I_p^*$. To shorten the
statements, we call singular fibers of all other types,
\ie, Kodaira's $\II$, $\III$, $\IV$ and $\II^*$, $\III^*$, $\IV^*$,
\emph{exceptional}. (These types are related to the exceptional
simple singularities/Dynkin diagrams $\bE_6$, $\bE_7$, $\bE_8$.)

Given two elliptic surfaces~$\XX_1$, $\XX_2$,
a fiberwise homeomorphism $\Gf\:\XX_1\to\XX_2$ is said to be
\emph{$2$-orientation preserving} (\emph{reversing}) if it
preserves (respectively, reverses) the complex orientation of the
bases and the fibers of the two elliptic fibrations.

\theorem\label{th.surface}
Two extremal
elliptic surfaces without exceptional fibers
are isomorphic if and only if
they are
related by a $2$-orientation preserving
fiberwise homeomorphism.
\endtheorem

Theorem~\ref{th.surface} is not proved separately,
as it
is an immediate consequence of
Theorem~\ref{th.hX} below:
the topological invariant distinguishing the surfaces is the
conjugacy class in $\bMG$ of the monodromy group
of the homological invariant~$\bhX$, see~\ref{s.j}.
In fact, we show that appropriate subgroups of $\bMG$
classify extremal elliptic surfaces without exceptional fibers,
both analytically and topologically.

Two extensions of Theorem~\ref{th.surface} to
somewhat wider classes
of surfaces are proved in
Subsections~\ref{s.ES} (see Remark~\ref{rem.type.IV})
and~\ref{s.rational}.


As
a by-product,
we obtain exponentially large collections
of non-home\-o\-mor\-phic elliptic surfaces sharing the same
combinatorial type of singular fibers.

\theorem\label{th.count}
For each integer $k\ge0$, there is a collection of\/ $\bT(k)$
extremal elliptic surfaces that share the same combinatorial type
of singular fibers, which is
\Dashes
\dash
$(k+2)\I_1\splus\I_{5k+4}^*$ if $k$ is even, or
\dash
$(k+2)\I_1\splus\I_{5k+4}$ if $k$ is odd,
\endDashes
but are not related by a $2$-orientation preserving
fiberwise homeomorphism.
\endtheorem

This theorem is proved in
Subsection~\ref{s.proof.count}, and
generalizations
are discussed in Subsection~\ref{s.more.surfaces}.
In fact, the
surfaces were constructed in~\cite{degt.kplets}, and
in~\cite{tripods} it was shown that they
share as well such topological invariants as
the transcendental lattice, see Example~\ref{ex.lattice},
and the
fundamental group of the ramification locus.

The proof of Theorems~\ref{th.surface} and~\ref{th.hX}
is based on an
explicit computation of the monodromy group $\Im\bhX$ of an
extremal
elliptic surface~$\XX$ in terms of its skeleton $\SkX$,
see~\ref{s.skeleton.def}.
In a sense,
we show that $\SkX$ \emph{is} $\Im\bhX$ (assuming that $\XX$ has
no type~$\II^*$ singular fibers).
As another consequence,
we obtain an algebraic
description of the reduced monodromy groups of such surfaces,
see Subsection~\ref{s.ES.group}, and a few results
(which may be known to the experts)
on the subgroups of the modular group~$\MG$; to me, the most
interesting seem Corollaries~\ref{cor.index}
and~\ref{free.product} describing the structure of subgroups
and Proposition~\ref{H.XY} characterizing monodromy groups of simple
\BM-factorizations.

\subsection{Real trigonal curves and real elliptic surfaces}
We consider a few other applications of the relation between
ribbon graphs and subgroups of~$\MG$, primarily to illustrate that
some classification problems are wilder than they may seem.

Recall that the \emph{Hirzebruch surface} is
the geometrically ruled surface $\Sigma_k\to\Cp1$, $k>0$, with an
exceptional section~$E$ of self-intersection $-k$. Up to
isomorphism, there is a unique real structure (\ie,
anti-holomorphic involution) $\conj\:\Sigma_k\to\Sigma_k$ with
nonempty real part $(\Sigma_k)_\R:=\Fix\conj$. A curve
$C\subset\Sigma_k$ is \emph{real} if it is invariant
under~$\conj$.
A \emph{trigonal
curve} is a curve $\CC\subset\Sigma_k$ disjoint from~$E$ and
intersecting each fiber of the ruling at three points. Such a
curve is \emph{generic} if all its singular fibers are of
type~$\I_1$ (simple tangency of the curve and a fiber
of the ruling). A generic
curve is necessarily nonsingular.

\theorem\label{th.curves}
For each integer $k\ge0$, there is a collection of
$\bT(k)$
generic real
trigonal curves $C_i\subset\Sigma_{2k+2}$
such that all real parts
$(C_i)_\R\subset(\Sigma_{2k+2})_\R$
are isotopic but the curves are not related by an equivariant
$2$-orientation preserving fiberwise auto-homeomorphism
of~$\Sigma_{2k+2}$ preserving the orientation of the real part
$\Rp1$ of the base of the ruling.
\endtheorem

Theorem~\ref{th.curves} is proved in
Subsection~\ref{s.proof.curves}, and a generalization
is discussed
in Subsection~\ref{s.ribbon.curves}.
The real part of each curve~$C_i$ in Theorem~\ref{th.curves}
consists of a `long' component~$L$
isotopic to $E_\R$ (see~\ref{s.real.part})
and
$(5k+4)$ ovals, all in the same connected component of
$(\Sigma_{2k+2})_\R\sminus(L\cup E_\R)$.

For each curve~$C_i$ as in Theorem~\ref{th.curves}, the double
covering $\XX_i\to\Sigma_{2k+2}$ ramified at $C_i\cup E$ is a real
Jacobian elliptic surface. Since the curves~$C_i$ are
distinguished by the braid monodromy,
one has the following corollary.

\corollary\label{cor.Lefschetz}
For each integer $k\ge0$, there
are two collections
of $\bT(k)$
real Jacobian elliptic surfaces $\XX_i\to\Cp1$
such that all real parts
$(\XX_i)_\R$ are fiberwise homeomorphic
but the surfaces are not related by an equivariant
$2$-orientation preserving fiberwise homeomorphism
of~$\Sigma_{2k+2}$ preserving the orientation of the real part
$\Rp1$ of the base of the elliptic pencil.
\qed
\endcorollary

In other words,
each of the two
collections consists of $T(k)$ pairwise non-\hskip0pt isomorphic
directed
real Lefschetz fibrations of genus~$1$ in the sense
of~\cite{Nermin}. The real parts $(\XX_i)_\R$ can be described in
terms of the \emph{necklace diagrams}, see~\cite{Nermin}: they
are chains of $(5k+4)$
copies of the same stone, which is either
\stone{\stcirc} or \stone{\square}.

\subsection{Contents of the paper}
In Section~\ref{S.ES}, we introduce the basic objects and prove
principal technical results relating extremal elliptic surfaces,
$3$-regular ribbon graphs, and geometric subgroups of~$\MG$.
Section~\ref{S.generalization} deals with a few generalizations of
these results to wider classes of ribbon graphs/subgroups.
In Section~\ref{S.trees}, we introduce \emph{pseudo-trees}, which
are ribbon graphs constructed from oriented rooted binary trees.
It is this relation that is responsible for the exponential drowth
in most examples. Theorem~\ref{th.count} is proved here.
In Sections~\ref{S.BM} and~\ref{S.real}, we prove the results
concerning, respectively, simple \BM-factorizations and real
trigonal curves.
Finally, in Section~\ref{S.TL} we introduce the notion of
transcendental lattice of a \BM-factorization and consider a few
examples.


\section{Elliptic surfaces\label{S.ES}}

In this section, we introduce basic notions and prove principal
technical results: Corollary~\ref{cor.Sk=H} and
Theorem~\ref{th.val}, establishing a connection between
$3$-regular ribbon graphs and geometric subgroups of~$\MG$, and
Theorems~\ref{th.jX} and~\ref{th.hX}, relating extremal elliptic
surfaces, their skeletons, and monodromy groups.

\subsection{The modular group}\label{s.MG}
Let
$\CH=\Z a\oplus\Z b$ be a rank~$2$ free abelian group
with the skew-symmetric bilinear form
$\bigwedge^2\CH\to\Z$ given by $a\cdot b=1$. We fix the
notation~$\CH$, $a$, $b$ throughout the paper
and \emph{define}
$\bMG:=\SL(2,\Z)$
as the group $\Sp\CH$
of symplectic auto-isometries of~$\CH$; it is
generated by the
isometries $\X,\Y\:\CH\to\CH$ given (in the basis $\{a,b\}$ above)
by the matrices
$$
\X=\bmatrix-1&1\\-1&0\endbmatrix,\qquad
\Y=\bmatrix0&-1\\1&\phantom{-}0\endbmatrix.
$$
One has $\X^3=\id$ and $\Y^2=-\id$. If $c=-a-b\in\CH$, then $\X$ acts
\via\
$$
(a,b)\overset \X\to\longmapsto(c,a)
 \overset \X\to\longmapsto(b,c)
 \overset \X\to\longmapsto(a,b).
$$
The \emph{modular group} $\MG:=\PSL(2,\Z)$ is
the
quotient $\bMG/\!\pm\id$. We retain the notation $\X$, $\Y$ for
the generators of~$\MG$. One has
$$
\MG=\<\X\,|\,\X^3=1\>\mathbin*\<\Y\,|\,\Y^2=1\>\cong\CG3\mathbin*\CG2.
$$

A subgroup $H\subset\MG$ is called \emph{geometric} if it is
torsion free and of finite index.
Since $\MG=\CG3\mathbin*\CG2$, the factors generated by~$\X$ and~$\Y$, a
subgroup $H\subset\MG$ is torsion free if and only if it is
disjoint from the conjugacy classes~$\cc[\X]$ and~$\cc[\Y]$, or,
equivalently, if both~$\X$ and~$\Y$ act freely on the quotient
$\MG/H$.

Similarly, a subgroup $\bH\subset\bMG$ is called
\emph{geometric} if it is torsion free and of finite index. A
subgroup $\bH\subset\bMG$ is torsion free if and only if
$-\id\notin\bH$ and the image of~$\bH$ in~$\MG$ is torsion free.

\subsection{Extremal elliptic surfaces}\label{s.surfaces}
In this subsection,
we remind a few well known facts concerning
Jacobian elliptic surfaces. The
principal
references are~\cite{FM} or the original paper~\cite{Kodaira}. For
more details concerning skeletons, we refer to~\cite{degt.kplets}.

A \emph{Jacobian elliptic surface} is a compact complex
surface $\XX$ equipped with an elliptic fibration $\pr\:\XX\to\BB$
(\ie, a fibration with all but finitely many fibers nonsingular
elliptic curves)
and a distinguished section $E\subset\XX$ of~$\pr$.
(From the existence of a section it follows that $\XX$ has no
multiple fibers.) Throughout the paper we assume that surfaces are
\emph{relatively minimal}, \ie, that fibers of
the elliptic pencil
contain no $(-1)$-curves.

\paragraph\label{s.Sigma}
Each nonsingular fiber of a Jacobian elliptic surface
$\pr\:\XX\to\BB$ is an abelian group, and the multiplication
by~$(-1)$ extends through the singular fibers of~$\XX$. The
quotient $\XX/\!\pm1$ blows down to a geometrically ruled surface
$\Sigma\to\BB$ over the same base~$\BB$, and the double covering
$\XX\to\Sigma$ is ramified over the exceptional section~$E$
of~$\Sigma$ and a certain \emph{trigonal curve} $C\subset\Sigma$,
\ie, a curve disjoint from~$E$ and intersecting each
generic fiber of the
ruling at three points.

\paragraph\label{s.j}
Denote by $\Bs\subset\BB$ the set of regular values of~$\pr$, and
define the \emph{\rom(functional\rom) $j$-invariant}
$\jX\:\BB\to\Cp1$ as the analytic
continuation of the function $\Bs\to\C^1$ sending each
nonsingular fiber of~$\pr$ to its classical $j$-invariant
(divided by~$12^3$).
The surface~$\XX$ is called \emph{isotrivial} if $\jX=\const$.

The monodromy
$\bhX\:\pi_1(\Bs,b)\to\bMG=\Sp H_1(\pr\1(b))$, $b\in\Bs$,
of the locally trivial fibration
$\pr\1(\Bs)\to\Bs$ is called the \emph{homological invariant}
of~$\XX$. Its reduction $\hX\:\pi_1(\Bs)\to\MG$ is called the
\emph{reduced monodromy}; it
is determined by
the $j$-invariant.
Together, $\jX$ and~$\bhX$ determine~$\XX$
up to isomorphism, and any pair
$(j,\smash\bhh)$ that agrees in the sense just described
gives rise to a unique isomorphism class of
Jacobian elliptic surfaces.

\paragraph\label{s.extremal}
According to~\cite{MNori},
a Jacobian elliptic surface~$\XX$ is extremal if and only if it
satisfies the following conditions:
\roster
\item\local{ext.1}
$\jX$ has no critical values other than~$0$, $1$, and~$\infty$;
\item\local{ext.2}
each point in $\jX\1(0)$ has ramification index at most~$3$,
and each point in $\jX\1(1)$ has ramification index at most~$2$;
\item\local{ext.3}
$\XX$ has no singular fibers of types~$\I_0^*$, $\II$,
$\III$, or~$\IV$.
\endroster

\paragraph\label{s.skeleton.def}
The
\emph{skeleton}
of a non-isotrivial elliptic surface $\pr\:\XX\to\BB$ is
the embedded bipartite graph
$\SkX:=\jX\1[0,1]\subset\BB$. The
pull-backs
of~$0$ and~$1$ are
called \black-- and \white-vertices of~$\SkX$, respectively.
(Thus, $\SkX$ is the \emph{dessin d'enfants} of~$\jX$ in the
sense of Grothendieck; however, we reserve the word `dessin' for
the more complicated graphs describing arbitrary, not necessarily
extremal, surfaces, \cf\ Subsection~\ref{s.dessins}.)
\emph{A priori}, $\jX$ may have critical values in
the open interval $(0,1)$, hence
the edges of~$\SkX$ may meet at points other than \black-- or
\white-vertices. However, by a small fiberwise equisingular
deformation of~$\XX$ the skeleton~$\SkX$ can be made
\emph{generic} in the sense that the edges of~$\SkX$ meet
only at \black-- or \white-vertices and the valency of each
\black-- (respectively, \white--) vertex is~$\le3$
(respectively,~$\le2$).

The skeleton~$\SkX$ of an extremal elliptic surface~$\XX$ is
always generic. In addition, each region of~$\SkX$
(\ie, component of $\BB\sminus\SkX$) is a topological disk; in
particular, $\SkX$ is connected.
Furthermore, each region contains a single critical point
of~$\jX$, the critical value being~$\infty$.
Thus, in this case $\SkX$ can be
regarded as an abstract ribbon graph: patching the cycles
of~$\SkX$ with disks, one recovers
the topological surface $\BB$ and ramified covering
$\jX\:\BB\to\Cp1$; then, the analytic structure on~$\BB$ is
given by the Riemann existence theorem.
It follows that the skeleton~$\SkX$
of an extremal elliptic surface~$\XX$ determines
its $j$-invariant $\jX\:\BB\to\Cp1$ (as an analytic function);
hence the pair
$(\SkX,\bhX)$ determines~$\XX$.

\paragraph\label{s.3-regular}
The exceptional singular fibers of an elliptic surface~$\XX$ are
in a one-to-one correspondence with the \black-vertices of~$\SkX$
of valency $\ne0\bmod3$ and its \white-vertices of valency
$\ne0\bmod2$. Hence, if $\XX$ is extremal and without exceptional
fibers, all \black-- and \white-vertices of~$\SkX$ are of
valency~$3$ and~$2$, respectively. Since $\SkX$ is a bipartite
graph, its \white-vertices can be ignored, assuming that such a
vertex is to be inserted at the middle of each edge
connecting two \black-vertices. Under this convention, the
skeleton of an extremal elliptic surface without exceptional
fibers is a $3$-regular ribbon graph.
As explained above, each region of~$\SkX$ is a disk containing a
single singular fiber of~$\XX$.
Hence
$\SkX$ is a strict
deformation retract of~$\Bs$, and the homological invariant can be
regarded as an anti-homomorphism $\bhX\:\pi_1(\SkX)\to\bMG$.
It is explained in~\cite{tripods} (see also
Remark~\ref{rem.orientation} below)
that $\bhX$ can be encoded in
terms of orientation of~$\SkX$.

\subsection{Skeletons: another point of view}\label{s.skeleton}
Following~\cite{tripods},
we start with redefining a $3$-regular ribbon graph
as a set of ends of its edges. However, in the further
exposition we will make no
distinction between a graph in the sense of
Definition~\ref{def.skeleton} below and its geometric
realization (defined in the obvious way).
We will also redefine a few notions related to graphs (like
connectedness, paths, \etc.); each time, it is immediately obvious
that the new notions are equivalent to their topological
counterparts defined in terms of geometric realizations.

\definition\label{def.skeleton}
A \emph{$3$-regular ribbon graph}
is a collection $\Sk=(\CE,\op,\nx)$, where
$\CE=\CESk$ is a finite set, $\op\:\CE\to\CE$ is a free involution, and
$\nx\:\CE\to\CE$ is a free automorphism of order three.
The orbits of~$\op$ are called the \emph{edges} of~$\Sk$,
the orbits of~$\nx$ are called its \emph{vertices}, and the orbits
of $\nx\1\op$ are called its \emph{faces} or \emph{regions}.
(Informally, $\op$ assigns to an end the other and of the same
edge, and $\nx$ assigns the next end at the same vertex with
respect to its cyclic order constituting the ribbon graph
structure.)

A \emph{pointed $3$-regular ribbon graph} is a pair $(\Sk,\ee)$,
where $\ee\in\CESk$.
\enddefinition

\Remark
Alternatively, one can consider $\CESk$ as the set of edges
of~$\Sk$ regarded as a bipartite ribbon graph,
see~\ref{s.3-regular}. Then the orbits of~$\op$ and~$\nx$
represent, respectively, the \white-- and \black-vertices
of~$\Sk$. Considering bipartite ribbon graph with the valency of
\black-- and \white-vertices equal to (respectively, dividing) two
given integers~$p$ and~$q$, one can extend, almost literally,
the material of
this and next subsections (respectively, the generalizations found
in Section~\ref{S.generalization}) to the subgroups of the group
$\<x,y\,|\,x^p=y^q=1\>$.
However, I do not know any interesting geometric applications of
this group.
\endRemark

\paragraph\label{s.MG-action}
Given a $3$-regular ribbon graph~$\Sk$,
the set~$\CESk$ admits a canonical left
$\MG$-action. To be precise, we define a homomorphism
$\MG\to\SG(\CESk)$
to the group $\SG(\CESk)$ of permutations of~$\CESk$
\via\ $\X\mapsto\nx\1$,
$\Y\mapsto\op$.
According to this convention, the vertices, edges, and regions
of~$\Sk$ are the orbits of~$\X$, $\Y$, and $\X\Y$, respectively.
The graph~$\Sk$ is \emph{connected} if and only
if the canonical $\MG$-action is transitive.
A connected $3$-regular ribbon graph is called a
\emph{$3$-skeleton}.

Given an element
$\ee\in\CESk$, we denote by $\Stab\ee\subset\MG$ its
stabilizer. Stabilizers of all points of a
$3$-skeleton form
a whole conjugacy class of subgroups of~$\MG$; it is denoted by
$\cc[\Stab\Sk]$ and is called the \emph{stabilizer} of~$\Sk$.

A \emph{morphism} of $3$-skeletons $\Sk'=(\CE',\op,\nx)$ and
$\Sk''=(\CE'',\op,\nx)$ is defined as a map $\CE'\to\CE''$
commuting with~$\op$ and~$\nx$. In other words, it is a morphism
of $\MG$-sets. A morphism of pointed $3$-skeletons $(\Sk',\ee')$ and
$(\Sk'',\ee'')$ is required, in addition, to take~$\ee'$ to~$\ee''$. The
group of automorphisms of a $3$-skeleton~$\Sk$ is denoted $\Aut\Sk$;
we regard it as a subgroup of the symmetric group $\SG(\CESk)$.

The following two statements, although crucial for the sequel,
are immediate consequences of the definitions.

\theorem\label{th.Sk=H}
The functors
$(\Sk,\ee)\mapsto\Stab\ee$,
$H\mapsto(\MG/H,H/H)$
establish an equivalence of the categories of
\Dashes
\dash
pointed $3$-skeletons and morphisms and
\dash
geometric subgroups $H\subset\MG$ and inclusions.
\qed
\endDashes
\par\removelastskip
\endtheorem

It follows that any morphism of $3$-skeletons is a
topological covering of their geometric realizations.

\corollary\label{cor.Sk=H}
The maps
$\Sk\mapsto\cc[\Stab\Sk]$,
$\cc[H]\mapsto\MG/H$
establish a canonical one-to-one correspondence
between the sets of
\Dashes
\dash
isomorphism classes of $3$-skeletons and
\dash
conjugacy classes of
geometric subgroups $H\subset\MG$.
\qed
\endDashes
\par\removelastskip
\endcorollary

If a $3$-skeleton~$\Sk$ is fixed, the isomorphism classes of pointed
$3$-skeletons $(\Sk,\ee)$ are naturally enumerated by the orbits of
$\Aut\Sk$. Hence one has the following corollary, concerning
properties of geometric subgroups.

\corollary\label{[H]}
The conjugacy class $\cc[H]$ of a geometric subgroup $H\subset\MG$
is in a one-to-one correspondence with the set of orbits of
$\Aut(\MG/H)$. Furthermore, there is an anti-isomorphism
$\Aut(\MG/H)=N(H)/H$,
where $N(H)$ is the normalizer of~$H$ \rom(acting on $\MG/H$ by
the right multiplication\rom).
\qed
\endcorollary

\Remark
Theorem~\ref{th.Sk=H}, as well as its
generalizations~\ref{th.Sk=H.3,1},~\ref{th.Sk=H.all}
below, relating subgroups of~$\MG$ and ribbon graphs
resemble the results of~\cite{Birch}. However, the two
constructions differ: in~\cite{Birch}, finite index subgroups of
the congruence subgroup $\MG(2)$ are encoded using bipartite ribbon
graphs with vertices of arbitrary valency.
Our approach is closer to that of~\cite{Bogomolov}, where the
modular $j$-function on a modular curve~$\BB$
(see~\cite{Shioda.modular} and Remark~\ref{rem.Shioda})
is described in terms of a special
triangulation of~$\BB$. Theorem~\ref{th.val} below and its
generalizations in Section~\ref{S.generalization} make the
geometric relation between ribbon graphs and subgroups of~$\MG$
even more transparent.
\endRemark

\subsection{Paths in a $3$-skeleton}\label{s.paths}
The treatment of paths found in~\cite{tripods} is
not quite satisfactory for our purposes; we choose a slightly
different approach here.

\definition\label{def.path}
A \emph{path} in a $3$-skeleton $\Sk=(\CE,\op,\nx)$ is a pair
$\Gg=(\ee,w)$, where $\ee\in\CESk$
and $w$ is a word in the alphabet $\{\op,\nx,\nx\1\}$.
The \emph{evaluation map} $\val$ sends a path $\Gg=(\ee,w)$ to
the element $\val\Gg\in\MG$ obtained by
replacing $\op\mapsto\Y$, $\nx^{\pm1}\mapsto\X^{\pm1}$ in~$w$
and multiplying in~$\MG$. The \emph{starting} and \emph{ending}
points of~$\Gg$ are, respectively, $\Gg_0:=\ee\in\CESk$ and
$\Gg_1:=(\val\Gg)\1\ee\in\CESk$.
A path~$\Gg$ is a \emph{loop} if $\Gg_0=\Gg_1$.
The \emph{product} of two paths $\Gg'=(\ee',w')$ and
$\Gg''=(\ee'',w'')$ is defined
whenever $\Gg''_0=\Gg'_1$; it is $\Gg'\cdot\Gg'':=(\ee',w'w'')$,
where $w'w''$ is the concatenation.
\enddefinition

\midinsert
\centerline{\picture{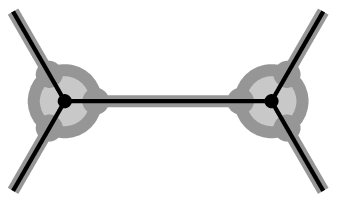}}
\figure\label{fig.path}
A $3$-skeleton~$\Sk$ (black), auxiliary graph~$\Ski$ (bold grey),
and space~$\Skii$ deformation equivalent to~$\Sk$
(bold and light grey)
\endfigure
\endinsert

\Remark\label{rem.path}
Intuitively, our definition of path represents the fact that, at
each point $\ee\in\CESk$, one can choose among three directions:
following the edge or walking around the vertex
preserving or reversing the cyclic order.
The inverse in the definition of~$\Gg_1$ is
due to the fact that the action of~$\MG$ is left rather than
right, hence the order of elements of~$w$ should be reversed.
(This is also one of the reasons why $\X$ is defined to act
\via\ $\nx\1$.) Strictly speaking, what is defined is a geometric
path (a chain of consecutive edges) in the auxiliary
graph~$\Ski$ obtained
from~$\Sk$ by shortening each edge and replacing each vertex with
a small circle (shown in bold grey lines in
Figure~\ref{fig.path}). The vertices of~$\Ski$ are in a natural
one-to-one correspondence with the elements of~$\CESk$.
When speaking about path homotopies,
fundamental groups, \etc., we replace~$\Ski$ with the topological
space~$\Skii$ obtained from~$\Ski$ by patching each circle with a
disk (light grey in the figure) and consider the homomorphisms
induced by the inclusion $\Ski\into\Skii$ and the strict
deformation retraction $\Skii\onto\Sk$.
\endRemark

The following two observations are also straightforward.

\lemma\label{lem.loop}
A path $\Gg$ is a loop if and only $\val\Gg\in\Stab\Gg_0$.
Conversely, given $\ee\in\CESk$, any element of\/ $\Stab\ee$ has the
form $\val\Gg$ for some loop $\Gg=(\ee,w)$.
\qed
\endlemma

\lemma\label{lem.mult}
Evaluation is multiplicative\rom:
$\val(\Gg_1\cdot\Gg_2)=\val\Gg_1\val\Gg_2$.
\qed
\endlemma

\theorem\label{th.val}
Given a pointed $3$-skeleton $(\Sk,\ee)$,
the evaluation map restricts to a well defined isomorphism
$\val\:\pi_1(\Sk,\ee)\to\Stab\ee$.
\endtheorem

\proof
Due to Lemmas~\ref{lem.loop} and~\ref{lem.mult}, it suffices
to show that $\val$ is well defined (\ie, it takes equal values on
homotopic loops) and $\Ker\val=\{1\}$.
Both statements follow from comparing the cancellations in
$\pi_1(\Sk,\ee)$ and in~$\MG$.

Since $\MG=\CG3\mathbin*\CG2$ is a free product, two words in
$\{\Y,\X,\X\1\}$ represent the same element of~$\MG$ if and only
if they are obtained from each other by a sequence of
cancellations of subwords of the form $\Y\Y$, $\X\X\1$, $\X\1\X$,
$\X\X\X$, or $\X\1\X\1\X\1$. The first three cancellations
constitute the
combinatorial definition of path homotopy in the
auxiliary graph~$\Ski$, see Remark~\ref{rem.path}: they correspond
to cancelling an edge immediately followed by its inverse. The
last two cancellations normally generate the kernel of the
inclusion homomorphism $\pi_1(\Ski,\ee)\to\pi_1(\Skii,\ee)$: they
correspond to contracting circles in $\Ski\subset\Skii$
to vertices of the original $3$-skeleton~$\Sk$.

An alternative proof of the fact that $\val$ is well defined is
given by Lemma~\ref{lem.h=val} below,
which provides an invariant geometric description of this map.
\endproof

\corollary\label{cor.geometric}
Any geometric subgroup $H\subset\MG$
\rom(respectively, any geometric subgroup
$\bH\subset\bMG$\rom)
has index divisible by six,
$[\MG:H]=6k$
\rom(respectively,
divisible by twelve,
$[\bMG:\bH]=12k$\rom)
and is
isomorphic to
a free group on $(k+1)$ generators.
\endcorollary

\proof
Let $\Sk=\MG/H$, see Theorem~\ref{th.Sk=H}. Then
$[\MG:H]=\ls|\CESk|$. On the other hand,
since $\Sk$ is a $3$-regular graph, one has $\ls|\CESk|=6k$ and $\Sk$ has
$2k$ vertices and $3k$ edges. Then $\chi(\Sk)=-k$ and $\pi_1(\Sk)$
is a free group on $(k+1)$ generators.

If $\bH\subset\bMG$ is a geometric subgroup, then
$\bH\not\ni-\id$ and the projection $\bH\to\MG$ is an
isomorphism onto its image, which is a geometric subgroup
of~$\MG$.
\endproof

\Remark\label{Farey.tree}
The universal covering
of a $3$-skeleton~$\Sk$ is a $3$-regular tree;
hence it is the Farey tree. The automorphism group $\Aut F$
of the Farey tree~$F$
can be identified with~$\MG$: it is generated by the rotations
about a vertex or the center of an edge. Thus, geometrically,
$\Sk=F/H$ for a finite index
subgroup $H\subset\Aut F$ acting freely on~$F$, and
Theorem~\ref{th.val} becomes a well known property of topological
coverings. If the action of~$H$ on~$F$ is not free, one needs to
consider the orbifold fundamental group $\piorb(F/H)$, see
Subsection~\ref{s.3,1} below.
If $[\MG:H]=\infty$, the quotient $F/H$ is
an infinite graph, see Subsections~\ref{s.infinite}
and~\ref{s.3,1.inf}.
\endRemark

\subsection{The homological invariant}\label{s.hX}
Fix a Jacobian elliptic surface $\pr\:\XX\to\BB$ without exceptional
fibers
and
let $\Sk=\SkX$ be the skeleton of~$\XX$. Assume that $\Sk$ is
generic, hence $3$-regular.
Consider the double covering $\XX\to\Sigma$ ramified at $\CC\cup E$,
see~\ref{s.Sigma}.
Pick a vertex~$v$ of~$\Sk$, let $F_v$ be the
fiber of~$\XX$ over~$v$, and let $\barF_v$ be its projection
to~$\Sigma$. Then, $F_v$ is the double covering of~$\barF_v$
ramified at $\barF_v\cap(\CC\cup E)$ (the three black points in
Figure~\ref{fig.basis} and~$\infty$).

\midinsert
\centerline{\picture{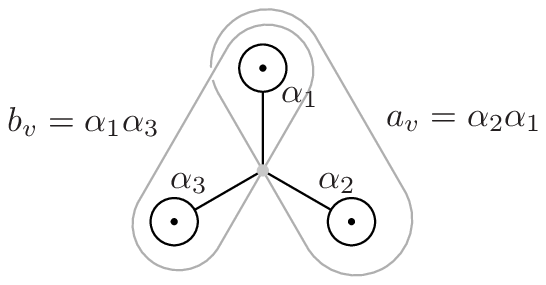}}
\figure
The basis in $H_1(F_v)$
\endfigure\label{fig.basis}
\endinsert

Recall that the three points of
intersection $\barF_v\cap\CC$ are in a
canonical one-to-one correspondence with the three ends
constituting~$v$, see~\cite{degt.kplets}.
Choose one of the ends
(a \emph{marking} at~$v$ in the terminology of~\cite{degt.kplets})
and let $\{\Ga_1,\Ga_2,\Ga_3\}$ be the
canonical basis for
the group $\pi_1(\barF_v\sminus(\CC\cup E))$ defined by
this end (see~\cite{degt.kplets} and Figure~\ref{fig.basis};
unlike~\cite{degt.kplets}, we take for the reference point the
zero section of~$\Sigma$, which is well defined in the presence of
a trigonal curve; this choice removes the ambiguity in the definition
of canonical basis).
Then $H_1(F_v)=\pi_1(F_v)$ is generated by the
lifts $a=\Ga_2\Ga_1$ and $b=\Ga_1\Ga_3$ (the two grey cycles in
the figure). To be precise, one needs to choose one of the two
pull-backs of the zero section
and take it
for the reference point for
$\pi_1(F_v)$
(the grey point at the center of the figure).
Thus, a choice of an end at~$v$ gives rise to an isometry
$H_1(F_v)=\CH$, which is canonical up to $\pm\id$.

Now, consider a copy~$F_\ee$ of~$F_v$ for each end $\ee\in v$ and
identify it with~$\CH$ using~$\ee$ as the marker.
(Alternatively, one can assume that a
separate fiber is chosen over each
vertex of the auxiliary graph~$\Ski$, see Remark~\ref{rem.path}.)
Under this
identification, the monodromy
$\bhh_\Gg\:H_1(F_{\Gg_0})\to H_1(F_{\Gg_1})$ of the locally trivial
fibration $\pr\1(\Sk)\to\Sk$ along a path~$\Gg$ in~$\Sk$ reduces to
a well defined element $\hh_\Gg\in\MG$.

\lemma\label{lem.h=val}
In the notation above, one has $\hh_\Gg=(\val\Gg)\1$.
\endlemma

\proof
Since both maps $\Gg\mapsto\hh_\Gg$ and
$\Gg\mapsto(\val\Gg)\1$ reverse products, it suffices to prove the
assertion for a path of length~$1$, \ie, for a single edge
of~$\Ski$.

Circumventing a vertex of the original skeleton~$\Sk$ in the
positive direction is the change of basis induced by a change of
the marker (rotation
through~$-2\pi/3$ about the center in Figure~\ref{fig.basis}); its
transition matrix is~$\X\1=(\val\,{\nx})\1$. Following an edge
of~$\Sk$ is a lift of the monodromy~$m_{1,1}$
in~\cite{degt.kplets}: during the monodromy, the black
ramification point
surrounded by~$\Ga_1$ crosses the segment connecting
the ramification points surrounded by~$\Ga_2$ and~$\Ga_3$;
modulo~$\pm\id$,
the corresponding
linear operator is given by $\Y=(\val\,{\op})\1$.
\endproof

Let $v$ be a vertex of~$\Sk$ and let $\ee\in v$. We will use the
notation $\pi_1(\Bs,\ee)$ for the group $\pi_1(\Bs,v)$, meaning that
the fiber~$F_v$ is identified with~$\CH$ using~$\ee$ as a marker.
Thus, we will speak about the reduced monodromy
$\hX\:\pi_1(\Bs,\ee)\to\MG$.

\theorem\label{th.jX}
Let $\XX$ be an extremal elliptic surface without exceptional
fibers, and let~$\ee$ be a
representative of a vertex of~$\SkX$.
Then the reduced monodromy $\hX\:\pi_1(\Bs,\ee)\to\MG$ takes values
in $\Stab\ee$, both maps in the diagram
$$
\pi_1(\SkX,\ee)\overset\inj_*\to\longrightarrow
 \pi_1(\Bs,\ee)\overset\hX\to\longrightarrow\Stab\ee\subset\MG
$$
are \rom(anti-\rom)isomorphisms, and the composed map is given
by $\Gg\mapsto(\val\Gg)\1$.
\endtheorem

\proof
Since $\SkX$
is a strict deformation retract
of $\Bs$, see~\ref{s.3-regular}, the
inclusion homomorphism $\inj_*\:\pi_1(\SkX)\to\pi_1(\Bs)$ is an
isomorphism. The rest follows from Lemma~\ref{lem.h=val} and
Theorem~\ref{th.val}.
\endproof

\theorem\label{th.hX}
The map $\XX\to\cc[\Im\bhX]$ establishes a
bijection
between the set of
isomorphism classes of extremal elliptic surfaces without
exceptional fibers and the set of
conjugacy classes of geometric subgroups of $\bMG$.
\endtheorem

\proof
It suffices to show that a subgroup $\bH\subset\bMG$ defines a
unique extremal elliptic surface. Since~$\bH$ is geometric, in
particular $-\id\notin\bH$,
the
projection $\bMG\to\MG$ induces an isomorphism of~$\bH$ to a
geometric subgroup $H\subset\MG$. The latter determines a skeleton
$\Sk\subset\BB$, hence a $j$-invariant
$j_X\:\BB\to\Cp1$ and corresponding
reduced monodromy $\hX\:\pi_1(\Bs)\to H$. Then, the
inverse isomorphism
$H\to\bH$ is merely a lift of~$\hX$ to a homological
invariant~$\bhX$;
together with~$j_X$, it defines a unique isomorphism class
of Jacobian elliptic surfaces, which are necessarily extremal due
to~\cite{MNori}, see~\ref{s.extremal}.
\endproof

Since the conjugacy class of the monodromy group of a fibration is
obviously invariant under fiberwise homeomorphisms,
Theorem~\ref{th.hX} implies Theorem~\ref{th.surface} in the
introduction.

\Remark
One can easily see that two extremal elliptic surfaces without
exceptional singular fibers are anti-isomorphic if and only if
their monodromy subgroups are conjugated by an element of
$\GL(2,\Z)\sminus\bMG$. (This conjugation results in a
homeomorphism of the skeletons reversing the cyclic order at each
vertex.) In other words, surfaces are anti-isomorphic if and only
if they are related by a $2$-orientation reversing homeomorphism.
\endRemark

\Remark\label{rem.Shioda}
The inverse map sending a
geometric subgroup $H\subset\bMG$ to an extremal
elliptic surface in
Theorem~\ref{th.hX} is equivalent to Shioda's
construction~\cite{Shioda.modular}
of modular elliptic surfaces, where the base~$\BB$ of the
elliptic fibration is obtained as the quotient
$\{z\in\C\,|\,\Im z>0\}/H$
and the $j$-invariant $\jX$ is the descent of the
modular $j$-invariant.
A generalization of the results of this section to arbitrary
finite index subgroups of~$\MG$
is considered in Subsections~\ref{s.3,1} and~\ref{s.ES}, see
Remark~\ref{rem.type.IV}; such
subgroups correspond to skeletons with monovalent \black-- and
\white-vertices allowed.
For a further generalization to arbitrary
subgroups, see Subsections~\ref{s.infinite} and~\ref{s.3,1.inf};
finitely generated subgroups can still be encoded by finite ribbon
graphs.
\endRemark

\Remark\label{rem.orientation}
In~\cite{tripods} it is shown that, for an extremal elliptic
surface~$\XX$ without exceptional singular fibers,
the homological
invariant~$\bhX$ admits a simple geometric description in terms of an
orientation of $\SkX$: one defines
the value $\bhX(\Gg)$ on a loop~$\Gg$
in~$\SkX$ to be
$\pm(\val\Gg)\1\in\bMG$,
depending on the parity of the number of edges
travelled by~$\Gg$ in the opposite direction. This correspondence
is not one-to-one, as distinct orientations may give rise to the
same homological invariant.
\endRemark

\section{Generalizations\label{S.generalization}}

In this section, we generalize some results of Section~\ref{S.ES}
to arbitrary subgroups of~$\MG$: finitely generated subgroups can
still be encoded by finite graphs.
Proofs are merely sketched, as
they repeat, almost literally, those in Section~\ref{S.ES}.
The material of this section is not used in the proofs of the
principal
results of the paper stated in the introduction.

\subsection{Infinite skeletons}\label{s.infinite}
In order to study subgroups of~$\MG$ of infinite index, we modify
Definition~\ref{def.skeleton} and define a \emph{generalized
$3$-regular ribbon graph} as a triple $\Sk=(\CESk,\op,\nx)$,
where $\CESk$ is a
set (not necessarily finite) and $\op$ and $\nx$ are free automorphisms
of~$\CESk$ of order~$2$ and~$3$, respectively. A \emph{generalized
$3$-skeleton} is a connected generalized $3$-regular ribbon graph.

All notions
introduced in Subsections~\ref{s.skeleton} and~\ref{s.paths} and
most statements proved there extend to the general case with
obvious changes. We restate Theorems~\ref{th.Sk=H}
and~\ref{th.val}.

\theorem\label{th.Sk=H.inf}
The functors
$(\Sk,\ee)\mapsto\Stab\ee$,
$H\mapsto(\MG/H,H/H)$
establish an equivalence of the categories of
\Dashes
\dash
pointed generalized $3$-skeletons and morphisms and
\dash
torsion free subgroups $H\subset\MG$ and inclusions.
\qed
\endDashes
\par\removelastskip
\endtheorem

\theorem\label{th.val.inf}
Given a pointed generalized $3$-skeleton $(\Sk,\ee)$,
the evaluation map restricts to a well defined isomorphism
$\val\:\pi_1(\Sk,\ee)\to\Stab\ee$.
\qed
\endtheorem

A generalized $3$-skeleton~$\Sk$ is called
\emph{almost contractible}
if
the group $\pi_1(\Sk)$
is finitely generated. Under Theorem~\ref{th.Sk=H.inf}, almost
contractible skeletons correspond to finitely generated torsion
free subgroups.

\proposition\label{almost.contractible}
There is a one-to-one correspondence between the sets of
\roster
\item
conjugacy classes of proper finitely generated torsion free subgroups
$H\subset\MG$,
\item
almost contractible $3$-skeletons with at least one cycle, and
\item\local{3.and.1}
connected finite ribbon graphs with all vertices of
valency~$3$ or~$1$
and such that distinct monovalent vertices are adjacent to distinct
trivalent vertices.
\endroster
Under this correspondence
$H\leftrightarrow\Sk\leftrightarrow\Skc$ one has
\rom(anti-\rom)isomorphisms
$N(H)/H=\Aut\Sk=\Aut\Skc$
and $H=\pi_1(\Sk)=\pi_1(\Skc)$\rom; in fact, $\Skc$ is
embedded to~$\Sk$ as an induced
subgraph and a strict deformation retract.
\endproposition

The finite ribbon graph~$\Skc$ corresponding to an almost
contractible $3$-skeleton $\Sk$ under
Proposition~\ref{almost.contractible}
is called the \emph{compact part} of~$\Sk$.
In the drawings, the monovalent vertices of~$\Skc$
(those that are to be
extended to `half' Farey trees) are represented by
triangles~\triang-,
\cf\ Figure~\ref{fig.trees2} in Subsection~\ref{s.two}.
The last condition
in~\iref{almost.contractible}{3.and.1} is the requirement that
$\Skc$ should admit no further contraction to a subgraph with all
vertices of valency~$3$ or~$1$. This condition makes $\Skc$
canonical.

\proof
Each almost contractible $3$-skeleton~$\Sk$
contains an induced subgraph~$\Sk'$ such that
$\Sk\sminus\Sk'$ is a forest: one can pick a finite collection of
loops
representing a basis for
$\pi_1(\Sk)$ and take for~$\Sk'$ the induced subgraph generated by
all vertices contained in at least one of the loops. (The
notation $\Sk\sminus\Sk'$ stands for the induced subgraph
generated by the vertices of~$\Sk$ that are not in $\Sk'$.)
The complement $\Sk\sminus\Sk'$ is a finite disjoint union of
\emph{infinite branches}, each infinite branch being a tree with
one bivalent vertex and all other vertices trivalent. Unless $\Sk$
is the Farey tree itself (corresponding to the trivial subgroup
of~$\MG$), each infinite branch is contained in a unique maximal
one. The maximal infinite branches are pairwise disjoint,
and contracting each such branch to its only bivalent vertex
produces the compact part~$\Skc$ as in the statement, the monovalent
vertices of~$\Skc$ corresponding to the maximal infinite branches
contracted. (The last condition
in~\iref{almost.contractible}{3.and.1} is due to the fact that, if
two monovalent vertices $u_1$, $u_2$ were adjacent to the same
vertex~$v$ then, together with~$v$, the two infinite branches
represented by~$u_1$ and~$u_2$ would form a larger infinite
branch.)

Since the construction is canonical, any automorphism of~$\Sk$
preserves~$\Skc$ and hence restricts to an automorphism of~$\Skc$.
Conversely, any automorphism of~$\Skc$ extends to a unique
automorphism of~$\Sk$: the uniqueness is due to the fact that
\emph{ribbon} graphs are considered; once an automorphism of such
a graph fixes a vertex~$v$ and an edge adjacent to~$v$, it is the
identity.
\endproof

\subsection{Skeletons with monovalent vertices}\label{s.3,1}
As another generalization, we lift the requirement that $\op$
and~$\nx$ should be free
and define a \emph{$(3,1)$-ribbon graph} as a triple
$\Sk=(\CESk,\op,\nx)$,
where $\CESk$ is a finite set and $\op$ and $\nx$ are automorphisms
of~$\CESk$ of order~$2$ and~$3$, respectively. A
\emph{$(3,1)$-skeleton} is a connected $(3,1)$-ribbon graph. Thus,
a $(3,1)$-skeleton is allowed to have monovalent \black-vertices
(which are the one element orbits of~$\nx$) and
`hanging edges' (one element
orbits of~$\op$); the latter are represented in the figures by
monovalent \white-vertices attached to these edges, \cf\
Figure~\ref{fig.path.3,1} below.

As above, all notions
introduced in Subsections~\ref{s.skeleton} and~\ref{s.paths}
extend to the case of $(3,1)$-skeletons.
Theorem~\ref{th.Sk=H} takes the following form.

\theorem\label{th.Sk=H.3,1}
The functors
$(\Sk,\ee)\mapsto\Stab\ee$,
$H\mapsto(\MG/H,H/H)$
establish an equivalence of the categories of
\Dashes
\dash
pointed $(3,1)$-skeletons and morphisms and
\dash
finite index subgroups $H\subset\MG$ and inclusions.
\qed
\endDashes
\par\removelastskip
\endtheorem

\paragraph\label{s.orbifold}
Denote by $D^2_1\cong D^2$, $D^2_2\cong\Rp1$, and~$D^2_3$ the
{\sl CW}-complexes obtained by
attaching a single $2$-cell $D^2$ to a circle~$S^1$ \via\ a map
$\partial D^2\to S^1$ of degree~$1$, $2$, or~$3$, respectively.
Given $\ee\in\CESk$,
define the \emph{orbifold fundamental group} $\piorb(\Sk,\ee)$
as the
fundamental group $\pi_1(\Skii,\ee)$, where
the space $\Skii$ is obtained
from~$\Sk$ by replacing a neighborhood of each trivalent
\black-vertex, monovalent \white-vertex, or monovalent
\black-vertex with a copy of~$D^2_1$, $D^2_2$, or~$D^2_3$,
respectively, \cf\ Figure~\ref{fig.path.3,1}.
(Note that $\piorb(\Sk,\ee)$ is indeed the orbifold fundamental
group, with the orbifold structure given by declaring each
monovalent \white-- or \black-vertex a ramification point of
ramification index~$2$ or~$3$, respectively. With this convention,
the universal covering of~$\Sk$ is again the Farey tree, \cf\
Remark~\ref{Farey.tree}.)
Contracting a maximal tree not containing a monovalent vertex, one
establishes a homotopy equivalence between~$\Skii$ and a wedge of
circles and copies of~$D^2_2$ and~$D^2_3$. Hence,
$\piorb(\Sk,\ee)$ is a free product
$$
\piorb(\Sk,\ee)=\circledast_{n_0}\Z\mathbin*
 \circledast_{n_2}\CG2\mathbin*\circledast_{n_3}\CG3,
\eqtag\label{eq.free.product}
$$
where $n_2$ and~$n_3$ are the numbers of monovalent \white-- and
\black-vertices, respectively, and $n_0=1-\chi(\Sk)=1-\chi(\Skii)$.
Observe that $\ls|\CESk|=6n_0+3n_2+4n_3-6$
(a simple combinatorial computation of the Euler characteristic).

\midinsert
\centerline{\picture{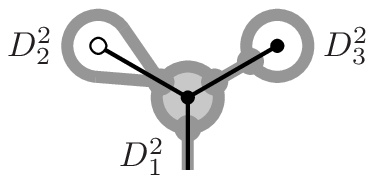}}
\figure\label{fig.path.3,1}
A $(3,1)$-skeleton~$\Sk$ (black), auxiliary graph~$\Ski$ (bold grey),
and space~$\Skii$ (bold and light grey)
\endfigure
\endinsert

Definition~\ref{def.path} of paths, loops, and the evaluation map
extends literally to the case of $(3,1)$-skeletons. Thus, we are
speaking about geometric paths in the auxiliary graph~$\Ski$
obtained by fattening the vertices of~$\Sk$ as shown in
Figure~\ref{fig.path.3,1}. (Note though that we disregard the
direction of a path
along the single edge replacing a \white-vertex and the
adjacent edge of~$\Sk$.) It is straightforward that $\piorb(\Sk)$
can be defined as the group of loops modulo an appropriate
equivalence relation. Next statement is proved similar to
Theorem~\ref{th.val}.

\theorem\label{th.val.3,1}
Given a pointed $(3,1)$-skeleton $(\Sk,\ee)$,
the evaluation map~$\val$ factors through a well defined isomorphism
$\val\:\piorb(\Sk,\ee)\to\Stab\ee$.
\qed
\endtheorem

\corollary\label{cor.index}
Any finite index subgroup $H\subset\MG$ is a free
product~\eqref{eq.free.product}, and one has
$[\MG:H]=6n_0+3n_2+4n_3-6$.
\qed
\endcorollary

\subsection{Extremal elliptic surfaces without type~$\II^*$
fibers}\label{s.ES}
Using the concept of $(3,1)$-skeleton introduced in the previous
section and
the description of the braid monodromy of the
ramification locus
found
in~\cite{degt.kplets} (the monodromy
$l_1(2)\mapsto\Y\X\1\Y$ and $l_1(3)\mapsto\Y$ for
monovalent \black-- and \white-vertices, respectively;
as in Subsection~\ref{s.hX}, the homomorphism $\BG3\to\MG$ is
given by~\eqref{eq.BG3toMG} below),
one
arrives at the following generalization of Theorem~\ref{th.jX}.

\theorem\label{th.jX.3,1}
Let $\XX$ be an extremal elliptic surface without type~$\II^*$
fibers, and let $\ee\in\CE$ be a
representative of a vertex of the skeleton~$\SkX$.
Then the reduced monodromy $\hX\:\pi_1(\Bs,\ee)\to\MG$ factors as
follows\rom:
$$
\pi_1(\Bs,\ee)\relbar\joinrel\onto
 \piorb(\SkX,\ee)\overset\cong\to\longrightarrow
 \Stab\ee\subset\MG
$$
where the rightmost anti-isomorphism is the map
$\Gg\mapsto(\val\Gg)\1$.
\qed
\endtheorem

\Remark\label{rem.Sk'}
In the presence of monovalent vertices, $\SkX$ is
no longer a
subspace of~$\Bs$. The first arrow in Theorem~\ref{th.jX.3,1} is
the composition of the homomorphisms induced by the strict
deformation retraction $\Bs\to\Sk'$ and the inclusion
$\Sk'\into\Skii$, where $\Sk'$ is obtained from~$\Ski$, see
Figure~\ref{fig.path.3,1}, by patching with disks the circles
surrounding the \emph{trivalent} \black-vertices only.
\endRemark

\corollary\label{cor.hX}
The map $\XX\to\cc[\Im\bhX]$ establishes a
bijection
between the set of
isomorphism classes of extremal elliptic surfaces without
type~$\II^*$ or~$\III^*$
fibers and the set of
conjugacy classes of finite index subgroups $\bH\subset\bMG$
such that $-\id\notin\bH$.
\endcorollary

\proof
Let~$X$ be a surface as in the statement, let
$\bH=\Im\bhX\subset\bMG$ (with respect to some base point
in~$\Bs$), and let $H=\Im\hX\subset\MG$ be the projection of~$\bH$
to~$\MG$. Under the assumptions, $\SkX$ has no \white-vertices and
hence $\piorb(\SkX)=H$ is a free product of copies of~$\Z$
and~$\CG3$ only. Furthermore, each order~$3$ generator of~$H$
represents the monodromy about a type~$\IV^*$
singular fiber of~$X$, see~\iref{s.extremal}{ext.3}, and hence
lifts to an order~$3$ element of~$\bH$. Thus, the projection
$\bH\to H$ admits a section and hence is an isomorphism. The rest
of the proof follows that of Theorem~\ref{th.hX}.
\endproof

\Remark\label{rem.type.IV}
Corollary~\ref{cor.hX} covers Shioda's
construction~\cite{Shioda.modular} to full extent and generalizes
Theorem~\ref{th.surface}
to surfaces with
type~$\IV^*$ fibers allowed.
Apparently, considering the homological invariant itself rather
than just its image, one can further generalize
Theorem~\ref{th.surface} to type~$\III^*$ singular
fibers. The special case of rational base is considered in
Theorem~\ref{th.1-1.P1} below.
\endRemark

\Remark
Surprisingly, type~$\II^*$ singular fibers do not fit into the
approach of this paper at all, as they are represented by bivalent
\black-vertices of the skeleton, \ie, orbits of~$\nx$ of length
two. Possibly, such skeletons can be treated as homogeneous spaces
of~$\bMG$ rather than~$\MG$, but the precise statements are not
quite clear at the moment.
An attempt of considering such more general skeletons
is made in~\cite{tripods}.
\endRemark

\subsection{The case of rational base}\label{s.rational}
In this subsection, we assume that the base~$\BB$ of an elliptic
fibration $\XX\to\BB$ is rational, $\BB\cong\Cp1$. In this case,
the homological invariant~$\bhX$ (lifting a given reduced
monodromy~$\hX$)
can be defined in terms of a
\emph{type specification} of~$\XX$, \ie, a choice of one of the two
possible types (whose local monodromies differ by $-\id$) of each
singular fiber. Moreover, the types of all but one singular fibers
can be chosen arbitrary, whereas the type of the remaining fiber
is determined by the requirement that the total multiplicity of
all singular fibers, which equals the topological Euler
characteristic $\chi(\XX)$, should be divisible by~$12$.
(The multiplicities of the two lifts of a given element of~$\MG$
differ by~$6$, \cf~\ref{s.degree}.)

If $\XX$ is extremal and has no type~$\II^*$ singular fibers, its
type specification can be described in terms of the reduced
monodromy group $H=\Im\hX$. Indeed, in view of
condition~\iref{s.extremal}{ext.3}, the types
of the exceptional fibers of~$\XX$ are fixed.
The non-exceptional singular fibers
are in a one-to-one correspondence with the regions
of~$\SkX$, equivalently, with the orbits of~$\X\Y$, equivalently,
with the $H$-conjugacy classes of maximal unipotent subgroups
of~$H$, and a type specification consists in assigning a lift
$\<\pm g\1(\X\Y)^ng\>\subset\bMG$ to each such conjugacy class
$\cc[\<g\1(\X\Y)^ng\>]_H$.

\theorem\label{th.1-1.P1}
Two extremal elliptic surfaces~$\XX_1$, $\XX_2$ over the rational
base $\BB=\Cp1$ and without type~$\II^*$ singular fibers are
isomorphic if and only if they are related by a $2$-orientation
preserving fiberwise homeomorphism.
\endtheorem

\proof
The `only if' part is obvious. For the `if' part, it suffices to
notice that a $2$-orientation preserving homeomorphism
$\XX_1\to\XX_2$ induces an orientation preserving homeomorphism
$\BB_1\to\BB_2$ taking punctures to punctures, commuting with the
homological invariants $\pi_1(\Bs_1)\to\MG\leftarrow\pi_1(\Bs_2)$
(and hence taking $H_1$ to~$H_2$) and preserving the type
specification (as distinct types of singular
elliptic fibers differ
topologically, for example by the local monodromy).
Hence, $\XX_1$ and~$\XX_2$ are isomorphic.
\endproof

\Remark\label{rem.ext.3}
The extremality condition
in Theorem~\ref{th.1-1.P1} can be relaxed by
replacing~\iref{s.extremal}{ext.3} by the requirement that the
surface should have no singular fibers of type $\I_0^*$, $\II^*$,
or~$\IV$. In this case, a type specification would also choose a
lift $\<\pm g\1\X g\>$ for each conjugacy class
$\cc[\<g\1\X g\>]_H$ of order~$3$ subgroups of~$H$ (monovalent
\black-vertices) and a lift $\<\pm g\1\Y g\>$ for each conjugacy
class $\cc[\<g\1\Y g\>]_H$ of order~$2$ subgroups of~$H$
(monovalent \white-vertices).
\endRemark

\Remark\label{rem.type}
The combinatorial type of singular fibers of an extremal (or more
general as in Remark~\ref{rem.ext.3}) elliptic surface~$\XX$ is
determined by its type specification and the following
combinatorial information about its skeleton~$\SkX$: the numbers
of monovalent \black-- and \white-vertices and the shapes of the
regions of~$\SkX$. Each monovalent \black-- (respectively,
\white--) vertex gives rise to a singular fiber of type $\II$
or~$\IV^*$ (respectively, $\III$ or~$\III^*$), and each $n$-gonal
region gives rise to a singular fiber of type~$\I_n$ or~$\I_n^*$.
There are large numbers of skeletons sharing these data; some
examples are considered in Subsections~\ref{s.proof.count},
\ref{s.more.surfaces}, and~\ref{s.sextics} below.
\endRemark

\subsection{The monodromy group of an elliptic surface}\label{s.ES.group}
For an elliptic surface~$\XX$, introduce the following fiber
counts:
\Dashes
\dash
$\nii$ is the number of fibers of type~$\II$ or~$\IV^*$;
\dash
$\niii$ is the number of fibers of type~$\III$ or~$\III^*$;
\dash
$\niv$ is the number of fibers of type~$\IV$ or~$\II^*$;
\dash
$t$ is the number of fibers of type~$\I_p^*$, $p\ge0$,
$\II^*$, $\III^*$, or~$\IV^*$.
\endDashes
Let, further, $\chi(\XX)$ be the topological Euler characteristic
of~$\XX$.

\theorem\label{th.Im.hX.max}
Let~$\XX$ be an extremal elliptic surface without type~$\II^*$
singular
fibers. Then the reduced monodromy group
$\Im\hX\subset\MG$ is a subgroup of index
$\chi(\XX)-6t-2\nii-3\niii$ isomorphic to the free product
$$
\circledast_{n}\Z\mathbin*
 \circledast_{\niii}\CG2\mathbin*\circledast_{\nii}\CG3,
$$
where $n=\frac16\chi(\XX)-t-\nii-\niii+1$.
\endtheorem

\proof
The statement follows from Theorem~\ref{th.jX.3,1},
Corollary~\ref{cor.index}, and the fact that
$$
\chi(\XX)=\ls|\CESk|+6t+2\nii+3\niii+4\niv,
\eqtag\label{eq.chi}
$$
where $\Sk=\SkX$.
(Here, we admit skeletons with bivalent \black-vertices as well.)
For the latter, observe that $\chi(\XX)$ equals the total
multiplicity of the singular fibers of~$\XX$. Exceptional singular
fibers are accounted for by the mono- and bivalent \black-vertices
and monovalent \white-vertices of~$\Sk$. Besides, there is one
fiber of type~$\I_p$ or~$\I_p^*$ inside each $p$-gonal region
of~$\Sk$. The sum of all indices~$p$ is the total number of
corners of all regions of~$\Sk$, \ie, $\ls|\CESk|$. Finally, each
$*$-type fiber increases the total multiplicity by~$6$.
\endproof

\theorem\label{th.Im.hX}
Let~$\XX$ be a non-isotrivial elliptic surface without
type~$\II^*$ or~$\IV$
singular
fibers. Then the index of the reduced monodromy group
$\Im\hX\subset\MG$ of~$\XX$ divides
$\chi(\XX)-6t-2\nii-3\niii$. In particular, it is finite.
\endtheorem

\proof
Let $\Sk$ be the skeleton of~$\XX$. After a fiberwise
equisingular deformation of~$\XX$,
not necessarily small, one can assume
that $\Sk$ is generic and connected. (For the modifications of
skeletons resulting in deformations of surfaces,
see~\cite{degt.kplets} or~\cite{DIK.elliptic}.)
Hence $\Sk$ is a $(3,1)$-skeleton. This time, each region of~$\Sk$
may contain several singular fibers of~$\XX$. Hence, instead of
Theorem~\ref{th.jX.3,1}, one has a diagram
$$
\pi_1(\Bs,\ee)\leftarrow\joinrel\relbar\joinrel\rhook
 \pi_1(\Sk',\ee)\relbar\joinrel\onto
 \piorb(\SkX,\ee)\overset\cong\to\longrightarrow
 \Stab\ee\subset\MG
$$
(where $\Sk'$ is the auxiliary space introduced in
Remark~\ref{rem.Sk'})
and an inclusion $\Stab\ee\subset\Im\hX$.
It remains to observe that $[\MG:\Stab\ee]=\ls|\CESk|$ and that
\eqref{eq.chi} holds for any non-isotrivial surface~$\XX$.
\endproof

\Remark
The reduced monodromy group $\Im\hX$
of an isotrivial elliptic surface~$\XX$ is
either trivial or conjugate to the subgroup generated by~$\X$
or~$\Y$. In particular, $[\MG:\Im\hX]=\infty$. At present, I do
not know whether the index of $\Im\hX$ is necessarily finite
if $\XX$ is a non-isotrivial surface
with type~$\II^*$ or~$\IV$ singular fibers.
\endRemark

\subsection{Further generalizations}\label{s.3,1.inf}
Combined, the constructions of Subsections~\ref{s.infinite}
and~\ref{s.3,1} give rise to the notion of
\emph{generalized\/ {\rm(\ie, possibly infinite)} $(3,1)$-skeleton}.
Theorems~\ref{th.Sk=H.inf} and~\ref{th.Sk=H.3,1} would combine to
the following statement.

\theorem\label{th.Sk=H.all}
The functors
$(\Sk,\ee)\mapsto\Stab\ee$,
$H\mapsto(\MG/H,H/H)$
establish an equivalence of the categories of
\Dashes
\dash
pointed generalized $(3,1)$-skeletons and morphisms and
\dash
subgroups $H\subset\MG$ and inclusions.
\qed
\endDashes
\par\removelastskip
\endtheorem

The orbifold fundamental group $\piorb(\Sk,\ee)$ of a generalized
$(3,1)$-skeleton~$\Sk$ is defined as in~\ref{s.orbifold},
and Theorem~\ref{th.val.3,1} extends
to this case
literally. Since $\Skii$ is still homotopy equivalent to a wedge
of circles and copies of~$D^2_2$ and~$D^2_3$, one obtains the
following
corollary.

\corollary\label{free.product}
Any subgroup of\/ $\MG$ is a free product \rom(possibly
infinite\rom) of copies of
cyclic groups~$\Z$, $\CG2$, and~$\CG3$.
\qed
\endcorollary

\paragraph\label{s.almost.contractible}
Under Theorem~\ref{th.Sk=H.all}, finitely generated subgroups
correspond to \emph{almost contractible $(3,1)$-skeletons}, which
are defined as those with finitely generated group $\piorb(\Sk)$.
Following the proof of
Proposition~\ref{almost.contractible}, one can
easily show that any almost contractible $(3,1)$-skeleton~$\Sk$
representing a
finitely generated subgroup $H\subset\MG$, $H\ne\{1\}$ (so that
$\Sk$ is not the Farey tree),
admits a strict deformation
retraction to a canonically defined finite induced subgraph
$\Skc\subset\Sk$, called the \emph{compact part} of~$\Sk$,
with the following properties:
\roster
\item
all vertices of $\Skc$ are of valency~$3$ or~$1$;
\item
the monovalent vertices of~$\Skc$ are divided into three types:
\white-, \black-, or \triang-\
(the latter representing
maximal infinite branches of~$\Sk$);
\item
distinct \triang-vertices are adjacent to distinct trivalent
vertices.
\endroster
Under this correspondence
$H\leftrightarrow\Sk\leftrightarrow\Skc$ one has
(anti-)isomorphisms
$N(H)/H=\Aut\Sk=\Aut\Skc$
and $H=\piorb(\Sk)=\piorb(\Skc)$, where
$\piorb(\Skc)$
is
defined similar to $\piorb(\Sk)$,
as the fundamental group of the space $(\Skc)^\bullet$ obtained
from~$\Skc$
by replacing each monovalent
\white-- or \black-vertex with a copy of~$D^2_2$ or~$D^2_3$,
respectively.

\section{Pseudo-trees\label{S.trees}}

Here, we introduce and count admissible trees and related ribbon
graphs, called pseudo-trees; they are the principal source of most
exponentially large examples stated in the introduction.

\subsection{Admissible trees and pseudo-trees}\label{s.trees}
An embedded tree $\tree\subset S^2$ is called \emph{admissible} if
all its vertices have valency~$3$ (\emph{nodes}) or~$1$
(\emph{leaves}).
Two such trees are called \emph{isomorphic} if they are related by
an orientation preserving auto-homeo\-mor\-phism of~$S^2$.
Each admissible tree~$\tree$ gives rise to its \emph{associated
$3$-skeleton} $\Sk_\tree$: one attaches a small loop to each leaf
of~$\tree$, see Figure~\ref{fig.tree}, left.
A $3$-skeleton obtained in this
way is called a \emph{pseudo-tree}. Clearly, each
pseudo-\allowbreak tree is a
skeleton of genus~$0$; two pseudo-trees $\Sk_{\tree'}$
and~$\Sk_{\tree''}$ are isomorphic as ribbon graphs if and only if
the trees~$\tree'$ and~$\tree''$ are isomorphic.

\midinsert
\centerline{\picture{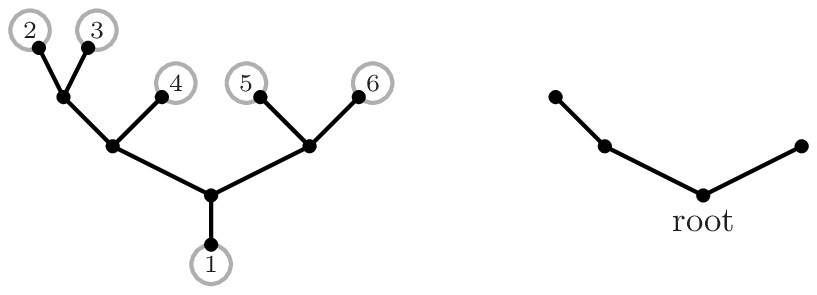}}
\figure
An admissible tree~$\tree$ (black) and
associated $3$-skeleton $\Sk_\tree$ (left);
the related binary tree (right)
\endfigure\label{fig.tree}
\endinsert

An admissible tree has a certain number $k\ge0$ of nodes and
$(k+2)$ leaves. The number of isomorphism classes of admissible
trees with $k$ nodes is denoted by~$\bT(k)$; it equals to the
number of isomorphism classes of pseudo-trees with $(2k+2)$
vertices.

\paragraph\label{s.marking}
A \emph{marking} of an admissible tree~$\tree$ is a choice of one of
its leaves~$v_1$. Given a marking, one can number all leaves
of~$\tree$ consecutively, starting from~$v_1$ and moving in the
clockwise direction (see Figure~\ref{fig.tree}, where the indices
of the leaves are shown inside the loops). Declaring the node
adjacent to~$v_1$ the root and removing all leaves, one obtains an
oriented rooted binary tree with $k$ vertices,
see Figure~\ref{fig.tree}, right;
conversely, an oriented rooted binary tree~$B$ gives rise to a unique
marked admissible tree: one attaches a leaf~$v_1$ at the root
of~$B$ and an extra leaf instead of each missing branch of~$B$. As
a consequence, the number of
isomorphism classes of
marked admissible trees with $k$
nodes is given by the Catalan number $C(k)$.

\paragraph\label{s.distance}
The \emph{vertex distance}~$m_i$ between two consecutive
leaves~$v_i$, $v_{i+1}$ of a marked admissible
tree~$\tree$ is the vertex length of the shortest left turn
path
in~$\tree$ from~$v_i$ to~$v_{i+1}$. For example,
in Figure~\ref{fig.tree} one has
$(m_1,m_2,m_3,m_4,m_5)=(5,3,4,5,3)$;
for another example, see Figure~\ref{fig.trees} in
Subsection~\ref{s.examples}.
The vertex distance between
two leaves~$v_i$, $v_j$, $j>i$, is defined to be
$\sum_{s=i}^{j-1}m_s$; it is the vertex length of the shortest
left turn path connecting~$v_i$ to~$v_j$ {\em in the associated
$3$-skeleton~$\Sk_\tree$}.

One can extend the sequence $(m_1,\ldots,m_{k+1})$ by appending the
vertex distance $m_{k+2}$ from~$v_{k+2}$ to $v_1$; then one has
$m_1+\ldots+m_{k+2}=5k+4$. Two marked trees are isomorphic if and
only if
their
sequences $(m_1,\ldots,m_{k+1})$ are
equal. Two unmarked trees are isomorphic if and only if the
corresponding extended sequences
$(m_1,\ldots,m_{k+1},m_{k+2})$ differ by a cyclic permutation.
Note that \emph{not} any sequence $(m_1,\ldots,m_{k+1})$ gives
rise to a marked admissible tree, see~\cite{tripods} for a
criterion.

\subsection{Counts}\label{s.counts}
As above, let
$\bT(k)$ be the number of isomorphism classes of pseudo-trees with
$(2k+2)$ vertices. Let, further, $\bT_i(k)$, $i\ge0$, be the
number of classes of pseudo-trees~$\Sk$ with $\ls|\Aut\Sk|=i$.

For a pseudo-tree~$\Sk$ with $(2k+2)$ vertices,
denote by $\gO_{\Sk}$ the orbit of
$\X\Y$
corresponding to the outer $(5k+4)$-gonal region
of~$\Sk$. The number of isomorphism classes of pointed $3$-skeletons
$(\Sk,\ee)$, where $\Sk$ is a pseudo-tree with $(2k+2)$ vertices
and $\ee\in\gO_{\Sk}$,
is denoted by~$\tT(k)$.

\lemma\label{counts}
For a pseudo-tree $\Sk=\Sk_\tree$ one has $\ls|\Aut\Sk|\le3$, \ie,
$\bT_i(k)=0$ for $i>3$.
The
numbers~
$T_1(k)$, $T_2(k)$, $T_3(k)$
are subject to the relations
$$
\gathered
\sum_{i=1}^3\frac{\bT_i(k)}i=\frac{C(k)}{k+2},\\\noalign{\vskip1\jot}
\bT_2(k)=\cases
 C(k'),&\text{if $k=2k'$},\\
 0,&\text{otherwise},
\endcases\qquad
\bT_3(k)=\cases
 C(k'),&\text{if $k=3k'+1$},\\
 0,&\text{otherwise}.
\endcases
\endgathered
$$
Furthermore, the group
$\Aut\Sk=\Aut\tree$ acts freely on the set of
leaves of
the original tree~$\tree$ and on the set~$\CESk$
of edge
ends of~$\Sk$.
\endlemma

\proof
Obviously, one has $\Aut\Sk_\tree=\Aut\tree$. Any combinatorial
automorphism of~$\tree$ is represented by a piecewise linear
auto-homeomorphism $\Gf\:\tree\to\tree$. Since $\tree$ is
contractible, $\Gf$ has a fixed point~$p$, which is necessarily
isolated (assuming that $\Gf\ne\id$, as an automorphism of
a connected ribbon graph fixing an edge is the identity). If $p$
is at the center of an edge of~$\tree$ (respectively, $p$ is a
vertex of~$\tree$), then $\Gf^2$ (respectively, $\Gf^3$) fixes a
whole edge of~$\tree$ and thus is the identity.

\midinsert
\centerline{\picture{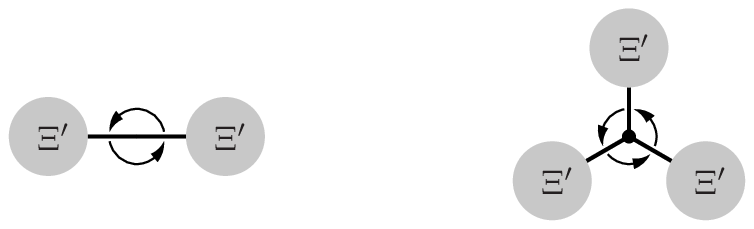}}
\figure
An automorphism of an admissible tree
\endfigure\label{fig.auto}
\endinsert

A tree~$\tree$ with an automorphism~$\Gf$ is shown in
Figure~\ref{fig.auto}. It is clear that such a tree admits no
automorphisms other than powers of~$\Gf$: the fixed point~$q$ of
such an automorphism would belong to one of the grey areas and the
vertices of~$\tree$ would be distributed unevenly about~$q$. Let
$k'$ be the number of nodes of the subtree~$\tree'$ shown in the
figure. In Figure~\ref{fig.auto}, left ($\ls|\Aut\tree|=2$), one
has $k=2k'$; in Figure~\ref{fig.auto}, right ($\ls|\Aut\tree|=3$),
one has $k=3k'+1$. In each case, the trees~$\tree$ admitting such
an automorphism~$\Gf$ can be parameterized by the marked
subtrees~$\tree'$, distinguished being the leaf extending towards
the fixed point of~$\Gf$. Their number is $C(k')$, which proves
the expressions for $\bT_2(k)$ and~$\bT_3(k)$.

It is also clear from Figure~\ref{fig.auto} that a non-trivial
automorphism does not fix a leaf of~$\tree$ or an edge end
of~$\Sk$. Then the first relation in the statement is the usual
orbit count: a tree~$\tree$ with $\ls|\Aut\tree|=i$ admits
$(k+2)/i$ essentially distinct markings, and the total number of
marked trees is $C(k)$.
\endproof

\corollary\label{cor.counts}
For each integer $k\ge0$, one has
$$
\bT(k)=\frac{C(k)}{k+2}+\frac{\bT_2(k)}2+\frac{2\bT_3(k)}3,\qquad
\tT(k)=\frac{5k+4}{k+2}C(k),
$$
where $\bT_2(k)$ and $\bT_3(k)$ are given by Lemma~\ref{counts}.
\endcorollary

\proof
Since $\bT_i(k)=0$ for $i>3$, the expression for
$\bT(k)=\bT_1(k)+\bT_2(k)+\bT_3(k)$ follows directly from
Lemma~\ref{counts}.

For each pseudo-tree~$\Sk$, one has $\ls|\gO_{\Sk}|=5k+4$ and
$\Aut\Sk$ acts freely on $\gO_{\Sk}$. Hence
$\tT(k)=(5k+4)\sum_{i=1}^3T_i(k)/i=(5k+4)C(k)/(k+2)$ due to
the first relation in Lemma~\ref{counts}.
\endproof

\subsection{Proof of Theorem~\ref{th.count}}\label{s.proof.count}
The surfaces in question were constructed in~\cite{degt.kplets}.
Each surface~$X$ corresponds to a pseudo-tree~$\Sk$
with $(2k+2)$ vertices, with
the type specification (see Subsection~\ref{s.rational} and
Remark~\ref{rem.type}) chosen so that the singular fiber
of~$X$
inside each monogonal region of~$\Sk$
should be of type~$\I_1$. The type of the
singular fiber inside the remaining $(5k+4)$-gonal region (the
outer region in Figure~\ref{fig.tree}, left) is then determined by
the parity of~$k$: it is of type $\I_{5k+4}$ if $k$ is odd or
$\I^*_{5k+4}$ if $k$ is even.

The $\bT(k)$ distinct pseudo-trees with $(2k+2)$ vertices give rise
to $\bT(k)$ pairwise non-isomorphic extremal elliptic surfaces;
Theorem~\ref{th.surface} implies that they are not related by a
$2$-orientation preserving fiberwise homeomorphism.
\qed

\subsection{Generalized pseudo-trees}\label{s.trees.a.c}
The construction of Subsection~\ref{s.trees} producing a $3$-skeleton
from a tree can be generalized. A function~$\ell$ defined on the
set of leaves of an admissible tree~$\tree$ and taking values
in~$\{0,\WHITE,\BLACK,\TRIANG\}$ is called \emph{admissible}
if no two leaves~$v_1$, $v_2$ with
$\ell(v_1)=\ell(v_2)=\TRIANG$
are adjacent to the same node. An
\emph{admissible pair} is a pair $(\tree,\ell)$, where $\tree$ is
an admissible tree and $\ell$ is an admissible function on the set
of leaves of~$\tree$. Each admissible pair $(\tree,\ell)$ gives
rise to an (almost contractible) $(3,1)$-skeleton
$\Sk_{(\tree,\ell)}$,
whose compact part $\Skc$ is obtained from~$\tree$
by attaching a small loop to each leaf~$v$ with $\ell(v)=0$ and
replacing each other leaf~$v$ with a monovalent vertex
of type $\ell(v)$, \cf\ Figures~\ref{fig.trees2}
and~\ref{fig.gtree} in Section~\ref{S.BM}.
Thus, one has
$\Sk_\tree=\Sk_{(\tree,0)}$. A generalized $(3,1)$-skeleton obtained in
this way is called a \emph{generalized pseudo-tree}.

Clearly, two generalized pseudo-trees $\Sk_{(\tree',\ell')}$ and
$\Sk_{(\tree'',\ell'')}$ are isomorphic if and only if so are
pairs $(\tree',\ell')$ and $(\tree'',\ell'')$, \ie, if there
exists an isomorphism $\Gf\:\tree'\to\tree''$ such that
$\ell'=\ell''\circ\Gf$.

For a generalized pseudo-tree $\Sk=\Sk_{(\tree,\ell)}$, we denote
by $\nstar(\Sk)$, $\STAR\in\{\WHITE,\BLACK,\TRIANG\}$,
the number of \emph{monovalent} \star-vertices of
the compact part~$\Skc$. Thus, $\nstar(\Sk)=\ls|\ell\1(\STAR)|$.

\proposition\label{H.XY}
Let $H\subset\MG$ be a proper finitely generated
subgroup. Then $H$ is
generated by $H\cap\cc[\X\Y]_\MG$ if and only if
$\MG/H$ is a generalized pseudo-tree without monovalent vertices
\rom(\ie, a skeleton $\Sk_{(\tree,\ell)}$ with $\ell$
taking values in $\{0,\TRIANG\}$\rom).
If this is the case,
$H$ admits a free basis
consisting of elements conjugate to $\X\Y$.
\endproposition

\proof
Let $\Sk=\MG/H$. It is an almost contractible $(3,1)$-skeleton,
see~\ref{s.almost.contractible}. Since $H$ is proper, $\Sk$ has a
well defined compact part~$\Skc$, which is \emph{not}
isomorphic to the skeleton
$\MGSk\WHITE$
representing~$\MG$ itself. Hence, each monogonal region
of~$\Sk$ (orbit of $\X\Y$ of length one)
is bounded by an edge with both ends attached to a
trivalent \black-vertex. (The only exceptional monogonal region is
the `outer' region in the skeleton $\MGSk\WHITE$
representing~$\MG$.) It
follows that the edge bounding a monogonal region
cannot belong to any subtree of~$\Skc$.

Let~$\tree$ be a maximal tree in~$\Skc$ not containing a monovalent
\white-- or \black-vertex. Contracting~$\tree$ establishes a homotopy
equivalence of the space $(\Skc)^\bullet$ computing
$\piorb(\Skc)=H$, see~\ref{s.almost.contractible}, to a wedge~$W$
of
circles and copies of~$D^2_2$ and~$D^2_3$. Each monogonal region
of~$\Sk$
produces a separate circle in~$W$, and the $H$-conjugacy classes
of loops represented by these circles constitute the intersection
$H\cap\cc[\X\Y]$. Thus, $H$ is generated by $H\cap\cc[\X\Y]$ if
and only if $W$ has no other circles or copies of $D^2_2$
or~$D^2_3$, \ie, $\Skc$ consists of several monogonal regions
attached to the (unique) maximal subtree $\tree\subset\Skc$.
\endproof

\Remark
Proposition~\ref{H.XY} gives a geometric characterization of the
proper
subgroups $H\subset\MG$ that can appear as the monodromy group of
a simple $\MG$-valued \BM-factorization, see
Definition~\ref{def.simple}.
Note that $\MG$ itself can also appear in this way (it is
generated by the images $\X\Y$ and $\X^2\Y\X\1$
of~$\Gs_1$ and~$\Gs_2$, respectively,
see~\eqref{eq.BG3toMG} below); it is the only monodromy
group that is not free.
\endRemark

\subsection{More examples of elliptic surfaces}\label{s.more.surfaces}
Let $\Sk=\Sk_{(\tree,\ell)}$ be a finite generalized pseudo-tree
(thus, we assume that $\ntriang(\Sk)=0$) obtained from
an admissible tree~$\tree$ with $k$ nodes. Let
$\nstar=\nstar(\Sk)$.
For the type specification (see Subsection~\ref{s.rational} and
Remark~\ref{rem.type}),
assign type~$\I_1$ to each monogonal region of~$\Sk$ and
types~$\IV^*$ and~$\III^*$ to the monovalent \black-- and
\white-vertices, respectively.
Then the fiber inside the remaining outer region of~$\Sk$ is of
type $\I_s$ if $k+\nblack+\nwhite$ is odd or $\I_s^*$ otherwise,
where $s=5k+4-\nblack-2\nwhite$.
(For even more examples, one could also vary the types~$\I_1$
or~$\I_1^*$ of the fibers in the monogonal regions, adjusting
the type of the
remaining
fiber
accordingly.)

The skeleton~$\Sk$ and the type specification described above
define an extremal elliptic surface~$\XX$ with the combinatorial
type of singular fibers
$$
(k+2-\nblack-\nwhite)\I_1\splus\nblack\IV^*\splus\nwhite\III^*
\splus\{\I_s\ \text{or}\ \I_s^*\}.
$$
The surfaces corresponding to non-isomorphic pairs $(\tree,\ell)$
are neither analytically isomorphic nor related by a
$2$-orientation preserving fiberwise homeomorphism, as they have
non-conjugate reduced monodromy groups.

\section{\BM-factorizations\label{S.BM}}

This section deals with \BM-factorizations. We prove
Theorems~\ref{th.strong} and~\ref{th.weak} and discuss a few
sporadic examples arising from generalized pseudo-trees and from
maximizing plane sextics.

\subsection{Preliminaries}\label{s.BM}
The \emph{braid group} $\BG3$ is the group
$$
\BG3=\<\Gs_1,\Gs_2\,|\,\Gs_1\Gs_2\Gs_1=\Gs_2\Gs_1\Gs_2\>=
 \<u,v\,|\,u^3=v^2\>,
$$
where $u=\Gs_2\Gs_1$ and $v=\Gs_2\Gs_1^2$. The center $Z(\BG3)$ is
the infinite cyclic group generated by $u^3=v^2$, and the quotient
$\BG3/Z(\BG3)$ is isomorphic to~$\MG$. In order to be consistent
with Subsection~\ref{s.hX}, we define the epimorphism
$\BG3\onto\bMG$ (and further to~$\MG$) \via
$$
\Gs_1\mapsto\X\Y,\quad
 \Gs_2\mapsto\X^2\Y\X\1.
\eqtag\label{eq.BG3toMG}
$$
(Then $u\mapsto-\X\1$ and $v\mapsto-\Y$.)

\paragraph\label{s.degree}
The abelianization $\BG3/[\BG3,\BG3]$ is the cyclic group~$\Z$.
The image of a braid
$\Gb\in\BG3$ in the abelianization $\BG3/[\BG3,\BG3]=\Z$ is called
its \emph{degree} $\deg\Gb$. (By convention, $\deg\Gs_1=1$.)
A braid $\Gb\in\BG3$ is uniquely recovered from its image
$\bar\Gb\in\MG$ and its degree $\deg\Gb$; the latter is determined
by~$\bar\Gb$ up to a multiple of~$6$. (The degree of an element
of~$\MG$ or~$\bMG$ is defined, respectively, modulo~$6$ or~$12$.)

\definition\label{def.simple}
A $\BG3$- (respectively, $\MG$- or $\bMG$-) valued
\BM-factorization $(\gm_i)$, $i=1,\ldots,r$, is called
\emph{simple} if each entry~$\gm_i$ belongs to the conjugacy class
$\cc[\Gs_1]$ (respectively, $\cc[\X\Y]_{\MG}$ or
$\cc[\X\Y]_{\bMG}$).
\enddefinition

\proposition\label{BG3=MG}
For each $r\ge1$,
the epimorphisms
$\BG3\onto\bMG\onto\MG$
establish
bijections between the sets of simple $\BG3$-, $\bMG$-, and
$\MG$-valued \BM-factorizations of length~$r$\rom; these
bijections preserve the weak/strong equivalence classes.
\endproposition

\proof
Each element
$x\in\cc[\X\Y]\subset\MG$ lifts to a unique element
$x'\in\cc[\Gs_1]\subset\BG3$
and to a unique element
$x''\in\cc[\X\Y]\subset\bMG$
(characterized by the requirement that $\deg x'=1$ and
$\deg x''=1\bmod12$),
establishing a one-to-one correspondence
between the sets of \BM-factorizations. The weak and strong
Hurwitz equivalences are preserved due to the fact that both
$\BG3\onto\MG$ and $\bMG\onto\MG$ are central extensions.
\endproof

\paragraph\label{s.length}
The advantage of considering the braid group~$\BG3$ rather than
the modular group $\MG$ is the fact that, in~$\BG3$, the
length~$r$ of a \BM-factorization of an element $\bminf\in\BG3$ is
uniquely determined by~$\bminf$: one has $r=\deg\bminf$. Hence,
for~$\BG3$,
the problem of uniqueness of a \BM-factorization
of a given element can be restated
in the language of \emph{factorization semigroup}, see~\cite{KK}
and~\cite{Orevkov.talk}.

\definition\label{def.BGM}
The \emph{factorization semigroup} is the semigroup~$\BGM{n}$
(with the group operation
denoted by~$\bmtimes$) generated by the elements
$\Gb\in\cc[\Gs_1]_{\BG{n}}$ subject to the \emph{Hurwitz relations}
$\Gb_1\bmtimes\Gb_2=\Gb_1\1\Gb_2\Gb_1\bmtimes\Gb_1=
\Gb_2\bmtimes\Gb_2\Gb_1\Gb_2\1$. The \emph{evaluation
anti-homomorphism} $v\:\BGM{n}\to\BG{n}$ is defined \via\
$v\:\Gb_1\bmtimes\Gb_2\bmtimes\ldots\bmtimes\Gb_r\mapsto\Gb_r\ldots\Gb_2\Gb_1$.
\enddefinition

\paragraph\label{s.B3toBn}
It is clear that an element $\gM\in\BGM{n}$ represents a strong
Hurwitz equivalence class of
simple $\BG{n}$-valued \BM-factorizations
(of length $\deg v(\gM)$)
and
the value $v(\gM)$ is merely the monodromy at infinity $\bmP\gM$.
Our Theorem~\ref{th.strong}
states that, for $n=3$, the evaluation map~$v$ is not injective;
moreover, the size of the pull-back
$v\1(\Gb)$, $\Gb\in\BG3$, may grow exponentially in the
degree $\deg\Gb$. Using the canonical inclusion $\BG3\into\BG{n}$,
one can easily conclude that the same assertion holds for any
integer $n\ge3$: the size of the pull-back
$v\1(\Gb)$, $\Gb\in\BG{n}$, may grow exponentially in the
degree $\deg\Gb$.

According to~\cite{Orevkov.talk},
the fact that $v$ is not injective implies that $\BGM{n}$ does
\emph{not} have the cancellation property, \ie,
an equality
$\Ga_1\bmtimes\Gb=\Ga_2\bmtimes\Gb$
or $\Gb\bmtimes\Ga_1=\Gb\bmtimes\Ga_1$
in~$\BGM{n}$
does not necessarily imply that $\Ga_1=\Ga_2$.

\subsection{Proof of Theorems~\ref{th.strong}
and~\ref{th.weak}}\label{s.proof.weak}
Consider a marked admissible tree $(\tree,v_1)$ with $k$ nodes and
$(k+2)$ leaves
and let $\Sk=\Sk_\tree$ be
the associated pseudo-tree, see
Subsection~\ref{s.trees}. Let $(m_1,\ldots,m_{k+1})$ be the
sequence of consecutive vertex distances, see~\ref{s.distance},
and consider the distances
$$
n_i=m_i+\ldots+m_{k+1},\ i=1,\ldots,k+1,\quad n_{k+2}=0
\eqtag\label{eq.n}
$$
from~$v_i$ to~$v_{k+2}$ in~$\Sk$.

Let $\ee\in\CESk$ be the edge end
at~$v_{k+2}$ that belongs to the original tree, see the grey dot
in Figure~\ref{fig.braids}, and consider the basis
$\{\Gg_1,\ldots,\Gg_{k+2}\}$ for $\pi_1(\Sk,\ee)$, where $\Gg_i$
is the class represented by the loop of~$\Sk$ attached at~$v_i$
which is
connected to~$\ee$ by the shortest left turn path in~$\Sk$ (the
grey loop in Figure~\ref{fig.braids}).

\midinsert
\centerline{\picture{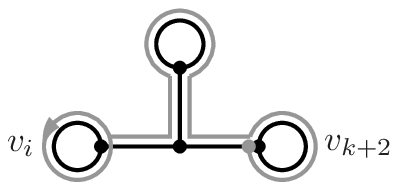}}
\figure
A loop $\Gg_i$ (grey)
\endfigure\label{fig.braids}
\endinsert

In terms of Definition~\ref{def.path}, the loop representing
a basis element~$\Gg_i$ is $(\ee,w_i)$, where
$$
w_i=(\nx\,\op)^{n_i}(\nx\,\op\,\nx\1\nx\1)(\op\,\nx\1)^{n_i}.
$$
The product $\Gg_1\ldots\Gg_{k+1}$ is homotopic to the boundary of
the outer $(5k+4)$-gonal region of~$\Sk$; after cancellation,
$\Gg_1\ldots\Gg_{k+2}\sim(\ee,(\nx\,\op)^{5k+4})$.

Define the $\MG$-valued \BM-factorization
$\gM=\gM(\Sk,\ee)=(\gm_1,\ldots,\gm_{k+2})$ \via
$$
\gm_i=(\val\Gg_i)\1=(\X\Y)^{n_i}(\X^2\Y\X\1)(\X\Y)^{-n_i}.
\eqtag\label{eq.MG}
$$
By construction, one has $\bmP\gM=(\X\Y)^{-5k-4}$, see
Lemma~\ref{lem.mult}, and
$\bmH\gM=\pi_1(\Sk,\ee)=\Stab\ee$, see
Theorem~\ref{th.hX}. Regarding each~$\gm_i$ in~\eqref{eq.MG} as an
element of~$\bMG$ and adjusting degree modulo~$12$, one
obtains $\bmP\gM=-(-\X\Y)^{-5k-4}\in\bMG$.

\Remark\label{rem.Artin}
Note
that the particular choice of a basis $\{\Gg_i\}$ used above is
not very important: by Artin's theorem~\cite{Artin},
any other basis $\{\Gg_i'\}$
with the property
that each $\Gg_i'$ is conjugate to some $\Gg_j$ and
$\Gg_1'\ldots\Gg_{k+2}'=\Gg_1\ldots\Gg_{k+2}$ is obtained from
$\{\Gg_i\}$ by a sequence of Hurwitz moves; hence the resulting
\BM-factorization~$\gM'$ would be strongly equivalent to~$\gM$.
\endRemark

Now, observe that $\ee$ belongs to the orbit $\gO_{\Sk}$
introduced in~\ref{s.counts}. Let $\ee'\in\gO_{\Sk}$ be another
element of this orbit, $\ee'=(\X\Y)^s\ee$, and consider the
\BM-factorization
$\gM'=\gM(\Sk,\ee'):=(\X\Y)^s\gM(\Sk,\ee)(\X\Y)^{-s}$.
Clearly, one has $\bmP{\gM'}=(\X\Y)^{-5k-4}$ and
$\bmH{\gM'}=\pi_1(\Sk,\ee')=\Stab\ee'$. As above, the strong
equivalence class of $\gM(\Sk,\ee')$ does not depend on the
particular choice of a basis for $\pi_1(\Sk,\ee')$; for this
reason, we omit the reference to the marking of the
original tree~$\tree$ in the notation.

Considering all $\tT(k)$ pairwise non-isomorphic pairs
$(\Sk,\ee)$, $\ee\in\gO_{\Sk}$, see
Subsection~\ref{s.counts} and
Corollary~\ref{cor.counts}, one obtains $\tT(k)$ distinct
\BM-factorizations $\gM(\Sk,\ee)$; they differ by the monodromy
groups $\bmH{\gm(\Sk,\ee)}=\Stab\ee$, see Theorem~\ref{th.Sk=H}.
Disregarding the base points~$\ee$, one arrives at $\bT(k)$ weak
equivalence classes, which differ by the conjugacy class
$\cc[\bmH{\gm(\Sk,\ee)}]=\cc[\Stab\Sk]$, see
Corollary~\ref{cor.Sk=H}.

The transcendental lattices and fundamental groups of the
\BM-factorizations constructed above are computed
in~\cite{tripods}; for the former, see Example~\ref{ex.lattice}.
\qed

\Remark
The \BM-factorizations~\eqref{eq.MG} represent the reduced
homological invariants of the extremal elliptic surfaces
constructed in Subsection~\ref{s.proof.count}.
\endRemark

\subsection{Examples}\label{s.examples}
Thus, the $\bT(k)$ weak equivalence classes of \BM-factorizations
given by Theorem~\ref{th.weak} are numbered by the isomorphism
classes of admissible trees with $k$ nodes. They are given
by~\eqref{eq.MG}, where the sequence $(n_1,\ldots,n_{k+2})$ is
obtained from the vertex distances $(m_1,\ldots,m_{k+1})$ of the
tree, see~\ref{s.distance}. The lifts to simple $\BG3$-valued
\BM-factorizations are
$$
\gm_i=\Gs_1^{n_i}\Gs_2\Gs_1^{-n_i},\ i=1,\ldots,k+2,\qquad
\bminf=(\Gs_1\Gs_2)^{3(k+1)}\Gs_1^{-5k-4}.
\eqtag\label{eq.BG}
$$
(For $\bminf$, we multiply $\Gs_1^{-5k-4}$ by a
power of the central element $(\Gs_1\Gs_2)^3$ in order to match
the degree.)

\example\label{ex.k=4}
The simplest example of non-equivalent \BM-factorizations given by
Theorem~\ref{th.weak} is obtained when $k=4$. The four admissible
trees with four nodes and their vertex distances are shown in
Figure~\ref{fig.trees}. The fact that the resulting
\BM-factorizations are not equivalent can be proved directly,
using \GAP~\cite{GAP}. Let~$\gM$ be one of the \BM-factorizations,
let $H=\bmH\gM$ be its monodromy group, and let~$N$ be the
normalizer of~$H$ in~$\MG$. Then, as Corollary~\ref{[H]} predicts,
the index $[N:H]$ equals~$1$, $2$, and~$3$ for the trees in
Figure~\ref{fig.trees}, left, middle, and right, respectively. In
particular, the four groups belong to at least three distinct
conjugacy classes. The two groups corresponding to the two trees in
the middle (which are related by an orientation reversing
diffeomorphism of the sphere) are conjugate in $\PGL(2,\Z)$ but
not in~$\MG$.
\endexample

\midinsert
\def\ccpic#1#2{$\vcenter{\halign{##\cr\picture{#1}\cr\hss$(#2)$\hss\cr}}$}
\centerline{\ccpic{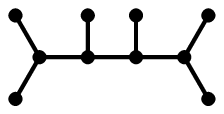}{3,4,4,4,3}\qquad\qquad
$\vcenter{\halign{#\cr\picture{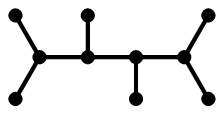}\cr\hss$(3,4,5,3,4)$\hss\cr
\noalign{\medskip}
\picture{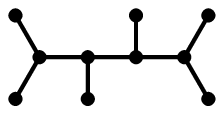}\cr\hss$(3,5,4,3,5)$\hss\cr}}$\qquad\qquad
\ccpic{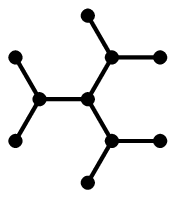}{3,5,3,5,3}}
\figure\label{fig.trees}
Admissible trees with four nodes
\endfigure
\endinsert

\example\label{ex.k=0}
The simplest example of weakly but not strongly equivalent
\BM-factorizations with the same monodromy at infinity
is given by Theorem~\ref{th.strong} with $k=0$.
The only admissible tree without nodes (two leaves connected by an
edge) gives rise to two \BM-factorizations:
$$
\gM'=(\Gs_1^2\Gs_2\Gs_1^{-2},\Gs_2),\qquad
\gM''=(\Gs_1\Gs_2\Gs_1\1,\Gs_1\1\Gs_2\Gs_1).
$$
Let $H',H''\subset\MG$ be their monodromy groups (reduced
to~$\MG$). Using \GAP~\cite{GAP}, one can see that
$[\MG:H']=[\MG:H'']=6$ whereas $[\MG:H'\cap H'']=24$. Hence
$H'\ne H''$.
\endexample

\subsection{Non-equivalent \BM-factorizations of length two}\label{s.two}
Consider the almost contractible
generalized pseudo-trees represented by the two
ribbon graphs shown in Figure~\ref{fig.trees2}. (Recall that each
\triang-vertex is to be extended to a maximal infinite branch,
which is a `half' of the Farey three, see
Subsection~\ref{s.infinite}.) They are obviously not isomorphic;
hence their stabilizers are not conjugate.

\midinsert
\centerline{\cpic{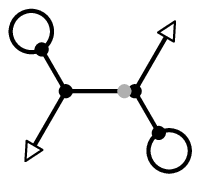}\qquad\qquad\qquad\cpic{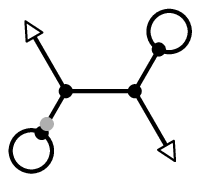}}
\figure\label{fig.trees2}
Almost contractible pseudo-trees with two loops
\endfigure
\endinsert

In each skeleton~$\Sk$, let $\ee\in\CESk$ be the edge end
represented by a grey dot in the figure, and pick a basis
$\{\Gg_1,\Gg_2\}$ for $\pi_1(\Sk,\ee)$ so that each $\Gg_i$,
$i=1,2$, is conjugate to the boundary of a monogonal region
of~$\Sk$ and $\Gg_1\Gg_2$ is homotopic to a circle encompassing
the compact part~$\Skc$ of~$\Sk$.
(The particular choice of bases is not important, see
Remark~\ref{rem.Artin}.)
Let
$\gM(\Sk)=((\val\Gg_1)\1,(\val\Gg_2)\1)$. For example, the bases
can be chosen so that
$$
\gather
\gM(\Sk_{\text{left}})=((\X\Y)(\X^2\Y\X\1)(\X\Y)\1,
 (\Y\X\Y)(\X^2\Y\X\1)(\Y\X\Y)\1),\\
\gM(\Sk_{\text{right}})=(\X^2\Y\X\1,
 (\Y\X\Y\X^2\Y)(\X^2\Y\X\1)(\Y\X\Y\X^2\Y)\1).
\endgather
$$
The $\BG3$-valued simple lifts of the two factorizations are
$$
\gM(\Sk_{\text{left}})=(\Gs_1\Gs_2\Gs_1\1,\Gs_2\Gs_1^3\Gs_2\Gs_1^{-3}\Gs_2\1),
\quad
\gM(\Sk_{\text{right}})=(\Gs_2,\Gb\Gs_2\Gb\1),
$$
where $\Gb=\Gs_2\Gs_1^2\Gs_2\1\Gs_1$.
One has
$$
\bmP{\gM(\Sk_{\text{left}})}=\bmP{\gM(\Sk_{\text{right}})}
 =\Y\X(\X\Y)^{-3}\Y\X(\X\Y)^{-3}
$$
(or, respectively, $\bminf=(\Gs_2\Gs_1^3\Gs_2\Gs_1\1)^2\in\BG3$).
On the other hand,
the monodromy groups $\cc[\bmH{\gM(\Sk)}]=\cc[\Stab\Sk]$
are not conjugate in~$\MG$ (although they are conjugate in
$\PGL(2,\Z)$);
hence the two \BM-factorizations are not weakly
equivalent.

\Remark
The two pseudo-trees differ by an orientation
reversing auto-homeomorphism of the sphere. This fact implies
that the
corresponding Hurwitz curves and Lefschetz fibrations are
anti-isomorphic. Hence, the two \BM-factorizations have isomorphic
fundamental groups and transcendental lattices,
see~\ref{s.invariants}.
\endRemark

\subsection{Non-simple \BM-factorizations}\label{s.more.BM}
Let $\Sk=\Sk_{(\tree,\ell)}$ be a generalized pseudo-tree obtained
from an admissible tree~$\tree$ with $k$ nodes, see
Subsection~\ref{s.trees.a.c}.
Denote $\nstar=\nstar(\Sk)$ for
$\STAR\in\{\BLACK,\WHITE,\TRIANG\}$.

Consider an embedding $\Skc\subset S^2$, patch each monogonal
region of~$\Skc$ with a disk, and let~$\BB$ be a regular
neighborhood of the result. Denote by~$\Bs$ the punctured disk
obtained from~$\BB$ by removing a point inside each monogonal
region of~$\Sk$ and all monovalent \black-- and \white-vertices
of~$\Sk$, see the shaded area in Figure~\ref{fig.gtree}. There is
an epimorphism $\Gr\:\pi_1(\Bs)\onto\piorb(\Sk)$, \cf\
Theorem~\ref{th.jX.3,1}.

\midinsert
\centerline{\picture{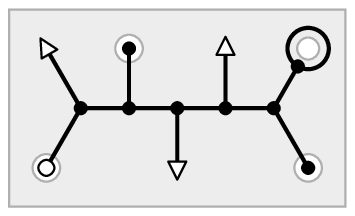}}
\figure\label{fig.gtree}
A generalized pseudo-tree~$\Skc$ and punctured disk~$\Bs$
\endfigure
\endinsert

Fix a point $b\in\partial\BB$ and pick a geometric basis
$\{\Gg_1,\ldots,\Gg_s\}$ for $\pi_1(\BB,b)$ such that
$\Gg_1\ldots\Gg_s=[\partial\BB]$. (The precise choice is not
important as different bases would produce weakly equivalent
\BM-factorizations, \cf\ Remark~\ref{rem.Artin}.) Define the
\BM-factorization $\gM(\Sk)=(\gm_1,\ldots,\gm_s)$ of length
$s=k+2-\ntriang$ \via\
$\gm_i=(\val\Gr(\Gg_i))\1$, $i=1,\ldots,s$.
It has $\nblack$ elements in $\cc[\X]$, $\nwhite$ elements
in~$\cc[\Y]$, and $k+2-\nblack-\nwhite-\ntriang$ elements
in~$\cc[\X\Y]$. Thus, $\gM$ is simple if and only if
$\nblack=\nwhite=0$.

If $\ntriang=0$, the conjugacy class of the monodromy at infinity
$\bmP{\gM(\Sk)}$
equals $\cc[(\X\Y)^{-n}]$, where $n=5k+4-\nblack-2\nwhite$, and
$\gM(\Sk)$ represents the reduced homological invariant of an
extremal elliptic surface constructed in
Subsection~\ref{s.more.surfaces}. In general, the monodromy at
infinity can be found as follows. Let $(m_1,\ldots,m_{\ntriang})$
be the sequence of vertex distances in~$\Skc$ between consecutive
\triang-vertices, each distance being the length of the shortest
left turn path connecting two \triang-vertices, with only
\black-vertices counted. (For example,
for the graph shown in Figure~\ref{fig.gtree},
starting from the upper left corner, one has
$(m_1,m_2,m_3)=(6,9,4)$; in Figure~\ref{fig.trees2}, for both
graphs one has $(m_1,m_2)=(5,5)$.) Then, the conjugacy class of
the monodromy at infinity $\bmP{\gM(\Sk)}$ is represented by the
\emph{right to left} product
$$
\prod_{i=1}^{\ntriang}(\X\Y)^{m_i-1}\X=
{}\ldots(\X\Y)^{m_2-1}\X\,\,(\X\Y)^{m_1-1}\X.
\eqtag\label{eq.minfty}
$$
Note that
$\sum_{i=1}^{\ntriang}m_i=5k+4-\nblack-2\nwhite-2\ntriang$.

\lemma\label{m-infty}
Given two generalized pseudo-trees~$\Sk'$, $\Sk''$, the
monodromies at infinity $\bmP{\gM(\Sk')}$ and $\bmP{\gM(\Sk'')}$
are conjugate in~$\MG$
if and only if the corresponding sequences $(m'_i)$ and
$(m''_j)$ differ by a cyclic permutation.
\endlemma

\proof
The `if' part is obvious. For the converse, observe that the
admissibility condition in Subsection~\ref{s.trees.a.c} implies
that each entry $m'_i$, $m''_j$ is at least~$2$. Then the
cyclic word~$w$
given by~\eqref{eq.minfty} admits no cancellations and the
distances $m_i$ can be recovered from the distances in~$w$ between
consecutive occurrences of~$\X^2$.
\endproof

\subsection{Maximizing plane sextics}\label{s.sextics}
We conclude this section with a few examples arising from
maximizing plane sextics.

Consider a plane sextic $\CC\subset\Cp2$ with simple singularities
only and with a distinguished type~$\bE$ singular point~$P$.
Let~$L_\infty$ be the (only) tangent to~$\CC$ at~$P$.
Assume that $L_\infty$ is not a component of~$\CC$ and let
$\CCa\subset\C^2=\Cp2\sminus L_\infty$ be the affine part of~$\CC$.
It is
a \emph{horizontal curve} in the sense of~\cite{Artal.braids} (or
\emph{Hurwitz curve} in the sense of~\cite{Kulikov}) of
degree~$3$ with respect to the pencil
$\CP=\{L_t\}$, $t\in\C^1$,
of lines through~$P$; in other words, the projection $\CCa\to\C^1$
defined by~$\CP$ is a proper map. Hence, using~$\CP$ and an
appropriately chosen section of the projection,
one can define the braid
monodromy $\bmC\:\pi_1(\Bs)\to\BG3$,
where $\Bs$ is the base~$\C^1$ of
the pencil with the singular fibers removed. Then, choosing a
geometric basis for $\pi_1(\Bs)$, one can represent~$\bmC$ by a
\BM-factorization~$\gMC$, which is well defined up to weak Hurwitz
equivalence.

The minimal resolution of singularities~$\XX$ of the double plane
ramified at a sextic~$\CC$ as above is a $K3$-surface, and the
pencil~$\CP$ lifts to an elliptic pencil $\XX\to\Cp1$
with a distinguished section. One can
easily show (see, \eg,~\cite{dessin-e7})
that $\XX$ is extremal if and only if $\CC$ is
\emph{maximizing}, \ie, if its total Milnor number takes its
maximal possible value~$19$. When this is the case, the
combinatorial type of singular fibers of~$\XX$ is determined by
the combinatorial type of singularities of~$\CC$ as follows:
\Dashes
\dash
the distinguished singular point~$P$ of type $\bE_6$, $\bE_7$,
or~$\bE_8$ gives rise to a singular fiber of type~$\I_6$, $\I_2^*$,
or~$\III^*$, respectively,
\dash
each other singular point gives rise to a singular fiber of the
following type: $\bA_p\mapsto\I_{p+1}$, $p\ge1$,
$\bD_q\mapsto\I_{q-4}^*$, $q\ge4$,
$\bE_6\mapsto\IV^*$, $\bE_7\mapsto\III^*$, $\bE_8\mapsto\II^*$,
\dash
a number of type~$\I_1$ fibers are added to make the total
multiplicity~$24$.
\endDashes
Furthermore, the $\bMG$-valued reduction of the braid
monodromy~$\bmC$ is the homological invariant~$\bhX$.

In~\cite{Artal.braids}, the authors construct a pair of reducible
maximizing sextics~$\CC_1$, $\CC_2$ with the set of singularities
$\bE_6\splus\bA_7\splus\bA_3\splus\bA_2\splus\bA_1$ and, using the
fact that both curves and all their singular fibers can be
chosen real, compute their \BM-factorizations $\gM_1$, $\gM_2$.
Then, reducing~$\gM_1$ and~$\gM_2$ to the finite group
$\SL(2,\CG{32})$ and using~\GAP~\cite{GAP}, they compute their
Hurwitz orbits and show that they are disjoint, concluding that
$\gM_1$ and~$\gM_2$ are not weakly equivalent and thus
distinguishing the curves.
(Both orbits are of length $15360$.)
In~\cite{dessin-e6}, the same pair of sextics is constructed using
trigonal curves or, equivalently,
extremal elliptic $K3$-surfaces; their
skeletons are as shown in Figure~\ref{fig.sextics}, with the
distinguished fiber~$L_\infty$ corresponding to the outer region.
Since the skeletons are obviously not isomorphic,
Theorem~\ref{th.hX} implies
that $\cc[\bmH{\gM_1}]\ne\cc[\bmH{\gM_2}]$.

\midinsert
\centerline{\cpic{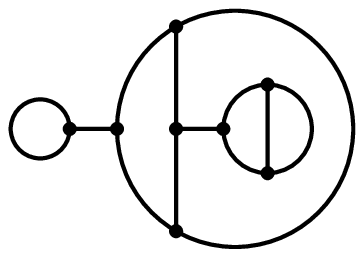}\qquad\qquad\qquad\cpic{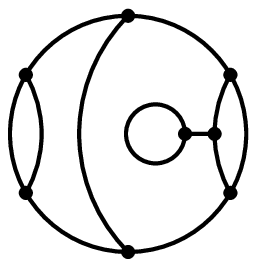}}
\figure\label{fig.sextics}
The set of singularities
$\bE_6\splus\bA_7\splus\bA_3\splus\bA_2\splus\bA_1$
\endfigure
\endinsert

\Remark
Strictly speaking, constructed in~\cite{dessin-e6} is
merely a pair of
not deformation equivalent sextics with the set of singularities
$\bE_6\splus\bA_7\splus\bA_3\splus\bA_2\splus\bA_1$.
However, it follows from~\cite{Shimada} that this set of
singularities is realized by exactly two equisingular deformation
families. Hence, the pairs found in~\cite{Artal.braids}
and~\cite{dessin-e6} coincide.
\endRemark

A number of other examples is found in~\cite{dessin-e8}
and~\cite{dessin-e6}. Listed in Table~\ref{tab.sextics} are all
sets of singularities realized by a pair~$\CC_1$, $\CC_2$ of
\emph{irreducible} maximizing plane sextics with a distinguished
type~$\bE$ singular point and with essentially different
skeletons. (More precisely, we ignore pairs of anti-isomorphic
curves.) For each such pair, Theorem~\ref{th.hX} implies
that the corresponding \BM-factorizations $\gM_1$, $\gM_2$ are not
weakly equivalent, as their monodromy groups are not conjugate.
For the sets of singularities marked with a~$^*$, the
corresponding \BM-factorizations differ by their transcendental
lattices, see Example~\ref{ex.sextics} below.

\midinsert
\table\label{tab.sextics}
Irreducible maximizing sextics with a type~$\bE$ singular point
\endtable
\def\*{\relax\llap{$^*$}}
\def\ba{\vtop\bgroup\halign\bgroup$##$\hss\cr}
\def\ea{\crcr\egroup\egroup}
\centerline{\ba
\* \bE_8\splus\bA_{10}\splus\bA_1\cr
 \bE_8\splus\bA_8\splus\bA_2\splus\bA_1\cr
\* \bE_8\splus\bA_6\splus\bA_4\splus\bA_1\cr
 \bE_8\splus\bA_5\splus\bA_4\splus\bA_2\cr
 (2\bE_6\splus\bA_5)\splus\bA_2\cr
 2\bE_6\splus\bA_4\splus\bA_3\cr
\ea\qquad\ba
 \bE_6\splus\bD_5\splus\bA_8\cr
\* \bE_6\splus\bD_5\splus\bA_6\splus\bA_2\cr
\* \bE_6\splus\bA_{10}\splus\bA_3\cr
\* \bE_6\splus\bA_{10}\splus\bA_2\splus\bA_1\cr
 \bE_6\splus\bA_9\splus\bA_4\cr
 \bE_6\splus\bA_8\splus\bA_4\splus\bA_1\cr
\ea\qquad\ba
 (\bE_6\splus\bA_8\splus\bA_2)\splus\bA_2\splus\bA_1\cr
\* \bE_6\splus\bA_7\splus\bA_4\splus\bA_2\cr
\* \bE_6\splus\bA_6\splus\bA_4\splus\bA_2\splus\bA_1\cr
\bE_6\splus\bA_5\splus2\bA_4\cr
(\bE_6\splus\bA_5\splus2\bA_2)\splus\bA_4\cr
\ea}
\endinsert


\Remark
It is worth mentioning that there also are three pairs~$\CC_1$,
$\CC_2$ of
irreducible maximizing sextics, those with the sets of
singularities
$$
\bE_7\splus\bE_6\splus\bA_4\splus\bA_2,\quad
\bE_7\splus\bA_{10}\splus\bA_2,\quad
\bE_7\splus\bA_6\splus\bA_4\splus\bA_2
$$
(the distinguished point~$P$ being that of type~$\bE_7$), such
that, within each pair, the curves
are not deformation equivalent but are represented by isomorphic
skeletons, hence have equivalent \BM-factorizations. It follows
that the affine parts $\CCa_1$, $\CCa_2$ are isotopic in the class
of Hurwitz curves, see~\cite{KK}. In fact, the curves constituting
each pair are related by a quadratic birational transformation
biholomorphic in the affine part $\Cp2\sminus L_\infty$.
\endRemark

\section{Real trigonal curves\label{S.real}}

Here, we give a brief introduction to theory of real trigonal
curves (see~\cite{DIK.elliptic} for more details),
prove Theorem~\ref{th.curves}, and consider a few
generalizations.

\subsection{Dessins}\label{s.dessins}
Recall that a \emph{real structure} on a complex analytic
variety~$\XX$ is an anti-holomorphic involution $\conj\:\XX\to\XX$.
A map, subvariety, \etc\. is called \emph{real} if it commutes
with/is preserved by~$\conj$.

For each Hirzebruch surface $\Sigma_k\to\BB\cong\Cp1$, $k\ge1$,
fix a (unique up to automorphism) real structure
$\conj\:\Sigma_k\to\Sigma_k$ with nonempty real part. Recall that
the ruling of~$\Sigma_k$ restricts to an $S^1$-fibration
$(\Sigma_k)_\R\to\BB_\R\cong\Rp1\cong S^1$ of the real parts, which
is orientable if and only if $k$ is even. The real part $E_\R$ of
the exceptional section $E\subset\Sigma_k$ is a section of this
fibration.

In what follows, we fix
an orientation of~$\BB_\R$ and denote by~$\BB_+$
the closure of the connected
component of $\BB\sminus\BB_\R$ whose complex orientation agrees
with the chosen orientation of the boundary
$\partial\BB_+=\BB_\R$.

\paragraph\label{s.def.dessin}
Given
a trigonal curve $\CC\subset\Sigma_k$, one can define
the \emph{$j$-invariant} $\jC\:\BB\to\Cp1$ by sending a
nonsingular fiber~$\barF$ to the $j$-invariant of the elliptic
curve~$F$ covering~$\barF$ and ramified at $\barF\cap(\CC\cup E)$.
(Here, the target is
the standard Riemann sphere $\C\cup\{\infty\}$.)
Following~\cite{Orevkov}
(see also~\cite{DIK.elliptic} for more details), define the
\emph{dessin} of~$\CC$ as the graph $\jC\1(\Rp1)\subset\BB$
with the following extra decoration:
\Dashes
\dash
the pull-backs of~$0$, $1$, and~$\infty$ are
\black--, \white--, and \cross-vertices, respectively;
\dash
the pull-backs of $[0,1]$, $[1,\infty]$, and $[-\infty,0]$ are
bold, dotted, and solid edges, respectively.
\endDashes
(Thus, the skeleton introduced in~\ref{s.skeleton.def} is obtained
from the dessin by removing all \cross-vertices and solid and
dotted edges.) The dessin of a real curve is invariant under the
complex conjugation in~$\BB$; for this reason, we only draw the
part contained in the closed disk $\BB_+$. Vertices and edges of
the dessin that belong to the boundary $\partial\BB_+$ are called
\emph{real}.

\Remark
Note that the $j$-invariant of a real curve may have real critical
values other than~$0$, $1$, or~$\infty$
not removable by
a small equivariant deformation. For this reason, a generic
symmetric dessin
may have non-removable \emph{monochrome} vertices in the
boundary~$\partial\BB_+$, \cf\ Figure~\ref{fig.dessin}.
\endRemark

According to~\cite{Orevkov} and~\cite{DIK.elliptic}, a dessin in
the \emph{topological} disk $\BB_+$ determines a real trigonal
curve~$\CC$, which is well defined up to equivariant fiberwise
deformation. (The converse is not true: a deformation of~$\CC$ may
result in a non-trivial modification of its dessin,
see~\cite{DIK.elliptic} for details. We do not use this fact
here.)

\paragraph\label{s.real.part}
From now on, we assume all curves nonsingular and generic, \ie, we
assume that all singular fibers are of Kodaira type~$\I_1$.

The real part $\CC_\R=\Fix\conj|_\CC$ of a real trigonal curve
$\CC\subset\Sigma_k$ consists of a \emph{long component}~$L$
isotopic to~$E_\R$ and a number of~\emph{ovals}, \ie, components
contractible in $(\Sigma_k)_\R$. The critical values of the
restriction $p\:\CC_\R\to\BB_\R$ of the ruling are the real
\cross-vertices of the dessin of~$\CC$. Pairs of such vertices
bound \emph{maximal dotted segments} in $\partial\BB_+$, each
segment containing a number of monochrome vertices and, possibly,
a number of real \white-vertices. The projection~$p$ is
three-to-one over the interior of each dotted segment, and it is
one-to-one outside the dotted segments. A maximal dotted segment
containing an even number of \white-vertices is the projection of
an oval, \cf\ Figure~\ref{fig.dessin}(a) and~(b); a segment
containing an odd number of \white-vertices is the projection of a
\emph{zigzag} in~$L$, \cf\ Figure~\ref{fig.dessin}(c).

The real \white-vertices of the dessin are the points where
$\CC_\R$ crosses the zero section of~$\Sigma_k$. It follows that,
if $k$ is even, two ovals of~$\CC_\R$ belong to the same
connected component
of the complement
$(\Sigma_k)_\R\sminus(L\cup E_\R)$ if and only if they are
separated by an even number of real \white-vertices.

\subsection{Proof of Theorem~\ref{th.curves}}\label{s.proof.curves}
To construct a curve~$\CC_i$ as in the statement,
consider one of the $\bT(k)$ pseudo-trees~$\Sk_i$ with $k$ nodes,
see Subsection~\ref{s.trees}, and extend it to a dessin as shown
in Figure~\ref{fig.dessin}(a) and~(b). More precisely,
embed~$\Sk_i$ to the sphere~$S^2$ (which is \emph{not} the base of
the elliptic pencil being constructed), patch each loop of~$\Sk_i$
with the disk bounded by this loop, and take for~$\BB_+$ a regular
neighborhood of the result in~$S^2$. Then, place a \white-vertex
at the center of each edge and a \cross-vertex at the center of
each disk bounded by a loop, connect all \black-- and
\white-vertices by appropriate edges
to the boundary~$\partial\BB_+$ in the radial manner, and
connect the resulting monochrome vertices in~$\partial\BB_+$
through \cross-vertices.

\midinsert
\centerline{\vbox{\halign{\hss#\hss&&\qquad\qquad\hss#\hss\cr
\cpic{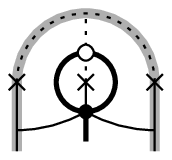}&\cpic{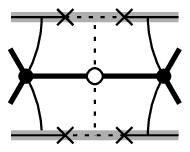}&\cpic{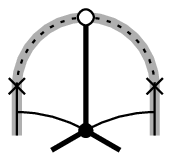}\cr
\noalign{\medskip}
(a)&(b)&(c)\cr}}}
\figure\label{fig.dessin}
Extending a skeleton~$\Sk$ (bold) to a dessin
\endfigure
\endinsert

Each loop of~$\Sk_i$ gives rise to an oval in~$\partial\BB_+$, see
Figure~\ref{fig.dessin}(a), and each edge of the original
tree~$\tree_i$ gives rise to two ovals, see
Figure~\ref{fig.dessin}(b).
Thus, we obtain the dessin of
a real trigonal curve in $\Sigma_{2k+2}$ with $(5k+4)$ ovals.
All curves obtained are topologically distinct: they differ by the
monodromy group $\cc[\Stab\Sk_i]$ of the monodromy
$\pi_1(\Bs_+)\to\MG$, where $\Bs_+$ is the interior of~$\BB_+$
with the inner \cross-vertices removed.
\qed

\subsection{A generalization: ribbon curves}\label{s.ribbon.curves}
The real trigonal curves constructed in
Subsection~\ref{s.proof.curves} are
\emph{ribbon curves} in the sense of~\cite{DIK.elliptic}. This
construction can be generalized.
Let $\Sk=\Sk_{(\tree,\ell)}$ be the generalized pseudo-tree obtained
from an admissible tree~$\tree$
with $k$ nodes
and a function~$\ell$ taking values in $\{0,\BLACK\}$,
see Subsection~\ref{s.trees.a.c}.
Let $z=\nblack(\Sk)$.
Extend~$\Sk$ do a dessin as shown in
Figure~\ref{fig.dessin}. The new element here is
Figure~\ref{fig.dessin}(c): the edge adjacent to a monovalent
\black-vertex~$v$ is extended towards $\partial\BB_+$ and $v$
is replaced with a \white-vertex (which is bivalent in the
complete dessin in~$\BB$),
giving rise to a zigzag rather than an
oval. The result is the dessin of a real trigonal curve
$\CC\subset\Sigma:=\Sigma_{2k+2-z}$ with $(5k+4-z)$ ovals
and $z$ zigzags.

\midinsert
\centerline{\picture{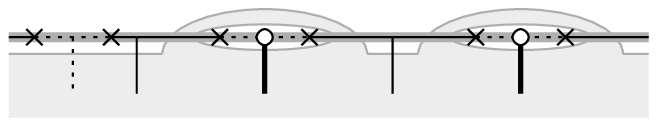}}
\figure\label{fig.zigzag}
The region $\Bs_+$
\endfigure
\endinsert

\paragraph\label{s.ribbon.topology}
To distinguish the curves topologically, consider the region
$\Bs_+$
obtained from the interior of~$\BB_+$ by adding small
regular neighborhoods of the zigzags and removing the zigzags
themselves and all inner \cross-vertices, see
Figure~\ref{fig.zigzag}. Since zigzags are clearly
distinguishable topologically, the monodromy $\pi_1(\Bs_+)\to\MG$
is a topological invariant of the curve. On the other hand,
at least topologically,
a pair of \cross-vertices constituting a zigzag can
collapse to a single type~$\II$ singular fiber; hence the
$\MG$-valued
monodromy about a whole zigzag equals that about a monovalent
\black-vertex. Thus, the image of the monodromy
$\pi_1(\Bs_+)\to\MG$ equals $\cc[\Stab\Sk]$, and distinct
skeletons produce non-isotopic curves.

\paragraph\label{s.ribbon.real.part}
Let $(m_1,\ldots,m_z)$
be the sequence
of vertex distances between consecutive
monovalent \black-vertices of~$\Sk$ (\cf\
Subsection~\ref{s.more.BM}; the monovalent
\black-vertices themselves are also included into the count,
so that each $m_i\ge3$).
Then the
topology of pair
$(\Sigma_\R,\CC_\R)$ is uniquely determined by the following two
properties:
\Dashes
\dash
$\CC_\R$ does not intersect the zero section except once inside
each zigzag;
\dash
the pair of zigzags of~$\CC_\R$
corresponding to a pair of consecutive
monovalent \black-vertices at a distance~$m$ is separated by
$(m-3)$ ovals.
\endDashes
Similar to Lemma~\ref{m-infty}, one can easily see that the
curves~$\CC'$, $\CC''$ obtained from two skeletons~$\Sk'$, $\Sk''$
as above have fiberwise isotopic real parts if and only if the
corresponding sequences $(m_i')$, $(m_j'')$ differ by a cyclic
permutation.

\paragraph\label{s.ribbon.surface}
If $z=\nblack(\Sk)$ is even, the double covering~$\XX$ of~$\Sigma$
ramified at~$\CC$ and~$E$ is a generic Jacobian real elliptic
surface. The surfaces obtained from distinct skeletons~$\Sk$ or
distinct (not related by an automorphism of~$\Sk$) lifts of the
real structure are neither deformation equivalent nor isomorphic
in the class of directed real Lefschetz fibrations, as they differ
by the homological invariants, \cf~\ref{s.ribbon.topology}. The
necklace diagram of~$\XX$, see~\cite{Nermin}, can be recovered
from the sequence $(m_1,\ldots,m_z)$ introduced
in~\ref{s.ribbon.real.part}: reading from~$m_z$ down to~$m_1$,
each pair $m_{2i}$, $m_{2i-1}$ gives rise to a copy of
\stone\stright, followed by $(m_{2i}-3)$ copies of \stone\stcirc,
a copy of \stone\stleft, and $(m_{2i-1}-3)$ copies of
\stone\square.
Two sequences produce isomorphic necklace diagrams if and only if
they differ by an \emph{even} cyclic permutation. (Thus, the lift
of the real structure is encoded in the choice of a marked
monovalent \black-vertex of~$\Sk$.)

\Remark
In the terminology of~\cite{DIK.elliptic}, the curves constructed
in this section are \emph{ribbon curves} with all blocks of
type~$\I_1$ or~$\II_3$. Conversely, any such curve~$\CC$ over the
rational base is obtained by the above construction, and the
ribbon curve structure of~$\CC$ is encoded by the original
skeleton~$\Sk$. It follows that both the fiberwise deformation
type and the fiberwise isotopy type of~$\CC$
determine
its ribbon curve structure. In~\cite{DIK.elliptic}, a similar
assertion is stated for ribbon curves with all blocks of
type~$\I_2$ or~$\II_3$.
\endRemark

\Remark
It is worth emphasizing that the analytic and topological
classifications of the curves constructed above coincide. This
fact substantiates the conjecture that real trigonal curves are
\emph{quasi-simple}, \ie, the fiberwise equisingular deformation
type of such a curve $\CC\subset\Sigma_k$ is determined by the
topological type of the quadruple $(\Sigma_k,\CC;\pr,\conj)$, where
$\pr\:\Sigma_k\to\Cp1$ is the ruling.
\endRemark

\section{The transcendental lattice\label{S.TL}}

In this section, we give a formal definition of a new invariant of
\BM-factorizations, which we call the \emph{transcendental
lattice}, and discuss a few open questions.

\subsection{The construction}\label{s.lattice}
Fix a commutative ring~$R$, two $R$-modules $\CL$, $\CV$,
and a skew-symmetric bilinear form
$\bigwedge^2\CL\to\CV$, $x\wedge y\mapsto x\cdot y$. (In case
$\CV$ has a
$2$-torsion, we assume, in addition, that $x\cdot x=0$ for all
$x\in\CL$.)
Fix, further, a symplectic (with respect to the chosen form)
representation $G\to\Sp\CL$.

\definition\label{def.lattice}
Given a $G$-valued \BM-factorization $\gM=(\gm_1,\ldots,\gm_r)$,
define the following objects:
\roster
\item
the $R$-module $\CL\gtimes\gM:=\bigoplus_{i=1}^r\CL$;
\item
the $R$-linear map $\chi\:\CL\gtimes\gM\to\CL$,
$\bigoplus_ix_i\mapsto\sum_i(\gm_i-1)x_i$;
\item
the $R$-quadratic map $q\:\CL\gtimes\gM\to\CV$,
$$
\bigoplus x_i\mapsto-\sum_{i=1}^rx_i\cdot\gm_ix_i+
 \sum_{1\le i<j\le r}(\gm_i-1)x_i\cdot(\gm_j-1)x_j.
$$
\endroster
(Here, $q$ is $R$-quadratic in the sense that $q(rx)=r^2q(x)$ for
all $x\in\CL\gtimes\gM$, $r\in R$ and
$(x,y)\mapsto q(x+y)-q(x)-q(y)\in\CV$ is a
$\CV$-valued bilinear form.)

Let $\CL_{\gM}=\Ker\chi$, and define
$\CL_{\gM}^\perp=\{x\in\CL_{\gM}\,|\,
 \text{$q(y+x)=q(y)$ for all $y\in\CL_{\gM}$}\}$. Then,
$\CL_{\gM}^\perp\subset\CL_{\gM}$ is an $R$-submodule and
the
quotient $\bmT\gM:=\CL_{\gM}/\CL_{\gM}^\perp$ inherits a
quadratic map $q\:\bmT\gM\to\CV$. It
is called the \emph{transcendental lattice} of~$\gM$
(defined by the representation $G\to\Sp\CL$).
\enddefinition

\lemma
One has $q(x+y)-q(x)-q(y)=\chi(x)\cdot\chi(y)\bmod2\CV$
for any pair $x,y\in\CL\gtimes\gM$.
\endlemma

\proof
The proof is a simple computation taking into account the fact
that each~$\gm_i$ is an isometry, so that
$\gm_ix_i\cdot\gm_iy_i+x_i\cdot y_i=2(x_i\cdot y_i)=0\bmod2\CV$.
\endproof

\corollary
If $\CV$ is free of $2$-torsion, the quadratic form
$q\:\CL_{\gM}\to\CV$ extends to a symmetric bilinear form
$\CL_{\gM}\otimes\CL_{\gM}\to\CV$.
\qed
\endcorollary

The symmetric bilinear extension of~$q$ is also denoted by~$q$.
Its kernel equals the submodule $\CL_{\gM}^\perp$
defined above,
and $q$ factors to a nondegenerate symmetric bilinear form
$q\:\bmT\gM\otimes\bmT\gm\to\CV$.
The pair $(\bmT\gM,q)$ is still called the
\emph{transcendental lattice} of~$\gM$.

\Remark
Assume that $\CL=H_1(F)$ for a punctured oriented
surface~$F$ and that the
map $G\to\Sp\CL$ is induced by a certain
representation of~$G$ in the
mapping class group of~$F$. In these settings, a weak Hurwitz
equivalence class of a G-valued \BM-factorization~$\gM$
of length~$r$
represents an $F$-bundle $X\to\Bs$ over a disk~$\Bs$ with $r$
punctures, see~\ref{s.geometry},
one has $\CL_{\gM}=H_2(X)$, and the symmetric bilinear form
$q\:\CL_{\gM}\otimes\CL_{\gM}\to\Z$ is given
by the intersection index,
$q\:x\otimes y\mapsto x\circ y$.
(Definition~\ref{def.lattice} is merely a
generalization of a simple algorithm computing $H_2(X)$ and the
self-intersections of $2$-cycles.)
The group $\CL\gtimes\gM$ can be interpreted as $H_2(X,F_b)$,
where $F_b$ is the fiber over a point $b\in\partial\Bs$, but the
form $q\:\CL\gtimes\gM\to\Z$ does not seem to have a geometric
meaning. Examples show that the associated bilinear form does not
need to be divisible by~$2$, see~\cite{tripods} or
Example~\ref{ex.lattice.ext} below.
\endRemark

\definition
A \emph{\rom(weak\rom) isomorphism} between two triples
$(\CM_1;\chi_1,q_1)$ and $(\CM_2;\chi_2,q_2)$, where $\CM_i$ is an
$R$-module, $\chi_i\:\CM_i\to\CL$ is an $R$-linear map, and
$q_i\:\CM_i\to\CV$ is an $R$-quadratic map,
is
an $R$-isomorphism $\Gf\:\CM_1\to\CM_2$ such that
$q_1=q_2\circ\Gf$ and $\chi_1=\chi_2\circ\Gf$ (respectively,
$\chi_1=g\circ\chi_2\circ\Gf$ for some $g\in G$).
\enddefinition

\proposition\label{tr.lattice}
The triples $(\CL\gtimes\gM;\chi,q)$ and
$(\CL\gtimes\gM';\chi',q')$ corresponding to two
strongly \rom(respectively, weakly\rom) equivalent
\BM-factorizations~$\gM$ and~$\gM'$ are isomorphic
\rom(respectively, weakly isomorphic\rom). In particular,
the transcendental lattice
$q\:\bmT\gM\to\CV$
is a weak equivalence invariant of~$\gM$.
\endproposition

\proof
If $\gM'$ is obtained from~$\gM$ by a global conjugation,
$\gm'_i=g\1\gm_ig$, $g\in G$, the weak isomorphism
$\CL\gtimes\gM'\to\CL\gtimes\gM$ is
$\Gf\:\bigoplus x_i'\mapsto\bigoplus gx_i'$; then
$\chi'=g\1\circ\chi\circ\Gf$.

Assume that $\gM'$ is obtained from $\gM$ by one
inverse Hurwitz move,
$$
\gm'_i=\gm_{i+1},\quad
\gm'_{i+1}=\gm_{i+1}\gm_i\gm_{i+1}\1,\quad
\gm'_j=\gm_j,\ j\ne i,i+1.
$$
Then the isomorphism
$\Gf\:\bigoplus x'_i\mapsto\bigoplus x_i$ is given by
$$
x_i=\gm_{i+1}\1x'_{i+1},\quad
x_{i+1}=x'_i+(\gm_i-1)\gm_{i+1}\1x'_{i+1},\quad
x_j=x'_j,\ j\ne i,i+1.
$$
It is straightforward that
$$
(\gm_i-1)x_i+(\gm_{i+1}-1)x_{i+1}=(\gm'_i-1)x'_i+(\gm'_{i+1}-1)x'_{i+1};
\eqtag\label{eq.sum}
$$
hence $\chi'=\chi\circ\Gf$. Furthermore, due to~\eqref{eq.sum},
the essentially different terms in the expressions for~$q$
and~$q'$ are
$$
-x_i\cdot\gm_ix_i-x_{i+1}\cdot\gm_{i+1}x_{i+1}+
 (\gm_i-1)x_i\cdot(\gm_{i+1}-1)x_{i+1}
$$
(and the corresponding primed terms).
Rewrite the latter sum in the form
$$
-x_i\cdot\gm_ix_i+[(\gm_i-1)x_i-x_{i+1}]\cdot(\gm_{i+1}-1)x_{i+1}
$$
(using $x_{i+1}\cdot x_{i+1}=0$) and observe that
$(\gm_i-1)x_i-x_{i+1}=-x'_i$ and
$$
(\gm_{i+1}-1)x_{i+1}=(\gm'_i-1)x'_i+(\gm'_{i+1}-1)x'_{i+1}
 -\gm_{i+1}\1(\gm'_{i+1}-1)x'_{i+1}.
$$
Multiplying out and using the
fact that $\cdot$ is skew-symmetric and $\gm_{i+1}=\gm'_i$
is an isometry, one obtains $q'=q\circ\Gf$.
\endproof

\subsection{Examples and open questions}\label{s.questions}
The transcendental lattice
$q\:\bmT\gM\to\CV$ is a relatively new invariant
(regarded as an invariant of a \BM-factorization) and I do not
know how powerful it is. In particular, I do not know if it can be
expressed in terms of other known invariants.

\problem
Is there a relation between $\bmT\gM$ and other known invariants,
for example $\cc[\bmH\gM]$ and $\cc[\bmP\gM]$?
\endproblem

Most known examples
of computation of~$\bmT\gM$
use the identity representation $\bMG=\Sp\CH$, see~\ref{s.MG}, and
deal with a \BM-factorization representing the
homological invariant of an extremal elliptic
surfaces~$\XX$. In this case, $\bmTL$ is indeed the transcendental
lattice of~$\XX$, \ie, the orthogonal complement
$\NS(\XX)^\perp\subset H_2(\XX)$,
with the form induced by the intersection index;
this relation explains the terminology, and it is the
computation
in~\cite{tripods} that
inspired Definition~\ref{def.lattice}.

\example\label{ex.sextics}
The $\bMG$-valued reductions of the
(non-simple)
\BM-factorizations arising from the pairs of
plane
sextics with the sets of singularities marked with a $^*$ in
Table~\ref{tab.sextics}, see Subsection~\ref{s.sextics}, differ by
their transcendental lattices. An easy way to prove this fact is
to compare the geometric classification of curves found
in~\cite{dessin-e8}, \cite{dessin-e6} and their arithmetic
classification found in~\cite{Shimada}. The same argument shows
that the other pairs in Table~\ref{tab.sextics} have isomorphic
transcendental lattices.
\endexample

\example\label{ex.lattice}
For each $k\ge0$, the simple $\bMG$-valued \BM-factorizations
given by Theorem~\ref{th.strong} have isomorphic transcendental
lattices, see~\cite{tripods}. If $k$ is even, one has
$\bmTL\cong\bD_k$ (with the usual convention $\bD_0=0$ and
$\bD_2=2\bA_1$); if $k=2s-1$ is odd, then $\bmTL$ is the
orthogonal complement
$(3\bv_1+\ldots+3\bv_s+\bv_{s+1}+\ldots+\bv_{2s-1})^\perp$
in the orthogonal direct sum
$\bigoplus_{i=1}^{2s-1}\Z\bv_i$, $\bv_1^2=1$.
\endexample

\paragraph
A \emph{coloring} of length~$r$
is a function
$\ell\:\{1,\ldots,r\}\to\{\pm1\}$.
Given a simple $\MG$-valued \BM-factorization
$\gM=(\gm_1,\ldots,\gm_r)$ and a coloring~$\ell$ of length~$r$,
define the lattice $\bmT{\gM,\ell}$ as the transcendental lattice
of the $\bMG$-valued lift of~$\gM$ obtained as follows: an entry
$\gm_i=g_i\1\X\Y g_i$, $i=1,\ldots,r$, $g_i\in\MG$,
lifts to $g_1\1\ell(i)\X\Y g_i\subset\bMG$.
Alternatively, this lift can be described as the one with the
eigenvalues of sign $\ell(i)$; in this form, the concept can be
extended to a wider class of \BM-factorizations, for example, to
those with unipotent entries, which arise from elliptic
surfaces/trigonal curves over the rational base and without
exceptional singular fibers. The following statement is immediate.

\proposition\label{prop.coloring}
Assume that two simple $\MG$-valued \BM-factorizations $\gM'$,
$\gM''$ of length~$r$ are weakly equivalent. Then there is a
permutation $\Gs\in\SG_r$ such that, for any coloring~$\ell$ of
length~$r$,
one has $\bmT{\gM',\ell}\cong\bmT{\gM'',\ell\circ\Gs}$.
\qed
\endproposition

\example\label{ex.lattice.ext}
In~\cite{tripods}, the lattices $\bmT{\gM,\ell}$ are computed for
all $\MG$-valued
\BM-factorizations given by Theorem~\ref{th.strong},
see~\eqref{eq.MG}, and all colorings~$\ell$ taking exactly one
value~$-1$. It turns out that the isomorphism class of
$\bmT{\gM,\ell}$ depends on $k$ only.
The corresponding quadratic forms $q\:\CH\gtimes\gM\to\Z$ are also
computed; in general, they do \emph{not} extend to integral
symmetric bilinear forms.
\endexample

\problem
Does Proposition~\ref{prop.coloring} distinguish the weak
equivalence classes given by Theorem~\ref{th.weak}?
\endproblem

\example
We conclude with the only known to me example of a direct
computation of the transcendental lattice using a representation
other than $\bMG=\Sp\CH$. In~\cite{ArimaShimada}, the authors
give an explicit construction of a pair of reducible
sextics (each splitting into an irreducible quintic~$Q$ and a
line~$L$) with the set of singularities $\bA_{10}\splus\bA_9$ and
compute their $\BG5$-valued braid monodromies with respect to the
pencil of lines through a generic point in~$L$. Then, following
more or less the lines of Definition~\ref{def.lattice} and using
the obvious representation $\BG5\to\Sp H_1(F)$, where $F$ is a
punctured surface of genus~$2$, they compute the transcendental
lattices and show that they are distinct (the latter fact being
predicted beforehand using theory of $K3$-surfaces). It is worth
mentioning that the two sextics are conjugate over $\Q(\!\sqrt5)$;
thus, $\bmTL$ is a topological, but not algebraic, invariant.
\endexample

\refstyle{C}

\enddocument